\documentclass[twoside,11pt,reqno]{amsart}
\usepackage{amsmath,amssymb,amscd,mathrsfs,epic}

\makeatletter

\@addtoreset{equation}{section}

\newtheorem{Proposition}{Proposition}[section]
\newtheorem{Lemma}[Proposition]{Lemma}
\newtheorem{Theorem}[Proposition]{Theorem}
\newtheorem{Corollary}[Proposition]{Corollary}

\newtheorem{Remark}[Proposition]{Remark}

\newtheorem{Example}[Proposition]{Example}

\newcommand{\infl}{\operatorname{infl}}

\newcommand{\Stab}{\operatorname{Stab}}
\newcommand{\lcm}{\operatorname{lcm}}

\newcommand{\ga}{\gamma}

\newcommand{\la}{\lambda}

\newcommand{\al}{\alpha}

\newcommand{\si}{\sigma}

\newcommand{\de}{\delta}

\newcommand{\eps}{\varepsilon}

\newcommand{\Ker}{\operatorname{Ker}}

\newcommand{\Irr}{{\mathrm {Irr}}}
\newcommand{\IBr}{{\mathrm {IBr}}}

\newcommand{\diag}{{\mathrm {diag}}}
\newcommand{\CC}{{\mathbb C}}

\newcommand{\OL}{{\mathcal O}}
\newcommand{\OLA}{{O_{\ell'}}}
\newcommand{\OLB}{{O_{\ell}}}
\newcommand{\SCL}{{\mathcal S}}

\newcommand{\ZZ}{{\mathbb Z}}
\newcommand{\NN}{{\mathbb N}}

\newcommand{\FF}{{\mathbb F}}
\newcommand{\FQ}{\FF_{q}}
\newcommand{\FQT}{\FF_{q}^{\times}}
\newcommand{\SSS}{{\sf S}}

\newcommand{\SC}{{L_{\CC}}}
\newcommand{\SF}{{S_{\FF}}}
\newcommand{\VC}{{V_{\CC}}}
\newcommand{\WN}{{W_{\CC}}}

\newcommand{\dar}{{\downarrow}}
\newcommand{\uar}{{\uparrow}}

\newcommand{\LF}{{\mathfrak{L}}}
\newcommand{\sfr}{{\mathfrak{s}}}
\newcommand{\tfr}{{\mathfrak{t}}}
\newcommand{\sfrs}{{\mathfrak{s}}^{\star}}

\newcommand{\tri}{\,{\trianglerighteq}\,}
\newcommand{\triop}{\,{\trianglelefteq}\,}

\renewcommand{\mod}{\bmod \,}

\begin{document}

\title[Representations of the special linear groups]
{{\bf Representations of finite special linear\\ groups in non-defining characteristic}}

\author{\sc Alexander S. Kleshchev}
\address
{Department of Mathematics\\ University of Oregon\\
Eugene\\ OR~97403, USA}
\email{klesh@uoregon.edu}

\author{\sc Pham Huu Tiep}
\address
{Department of Mathematics, University of Florida, 
Gainesville, FL~32611, USA}
\address{{\em Since August 2008:} Department of Mathematics, University of Arizona,
Tucson, AZ~85721, USA}
\email{tiep@math.arizona.edu}

\thanks{
\indent
Research supported by the NSF (grants DMS-0654147 and DMS-0600967). This work was completed when both authors visited MSRI, whom we thank for support.}


\maketitle

\section{Introduction}\label{SIntro}

\subsection{Background}
Let  $q=p^f$ be a power of a prime number $p$, and $SL_n(q)=SL_n(\FF_q)$ be the special 
linear group over the field $\FF_q$ with $q$ elements. Let $\FF$ be an algebraically closed 
field of characteristic $\ell\geq 0$. 
We are interested in parametrizing irreducible representations of $SL_n(q)$ over $\FF$.
If $\ell=p$, 
this is well known: a parametrization is given by the so-called $q$-restricted dominant weights, 
see for example \cite[Theorem~43]{Steinberg}.

Assume from now on that $\ell\neq p$. Then it seems natural to proceed through the general 
linear group  $GL_n(q)=GL_n(\FF_q)$. Indeed, a natural classification of irreducible 
$\FF GL_n(q)$-modules is available \cite{D1, D2, DJ1,J}, 
and one might try to study irreducible $\FF SL_n(q)$-modules 
by restricting from $GL_n(q)$ to $SL_n(q)$. A parametrization of complex irreducible characters 
of $SL_{n}(q)$ along these lines was obtained in \cite{L} (see also \cite{Bo}). 

Restricting irreducible modules from a finite group $G$ to a normal subgroup $S$ is the 
subject of Clifford theory. The first key general fact is that the restriction $V\dar_S$ 
of an irreducible $\FF G$-module $V$ is completely reducible. The situation is especially nice 
when $G/S$ is cyclic, e.g. $GL_n(q)/SL_n(q)\cong C_{q-1}$. In this case the restriction 
$V\dar_S$ is also multiplicity-free.
Clifford theory is most effective if $G/S$ is a cyclic group of order prime to $\ell$. So the 
problem of classifying irreducible $\FF SL_n(q)$-modules is easy if $\ell{\not{|}}(q-1)$. An 
additional consideration shows that the problem is still easy if $\ell{\not{|}}\gcd(n,q-1)$. 

However, in the difficult case where $\ell|\gcd(n,q-1)$ there are no general tools 
to describe the number of irreducible summands in the restriction $V\dar_{SL_n(q)}$ 
of an irreducible $\FF GL_n(q)$-module $V$. The main results of this paper are a 
precise description of this number and a corresponding parametrization of irreducible $\FF SL_n(q)$-modules. 

We note that even a classification of irreducible $\FF GL_n(q)$-modules which are 
{\em irreducible} on restriction to $SL_n(q)$, which is of course a special case of 
our problem, seems to have been unknown. This question is important for Aschbacher-Scott  
program \cite{A, Sc} on classifying maximal subgroups in finite classical groups 
and has been our original 
motivation, see \cite{KTNew}. Although representations of $SL_{n}(q)$ in non-defining 
characteristic were studied in \cite{Gru} using a different approach, it is not clear how 
to use \cite{Gru} to describe irreducible restrictions from $GL_n(q)$ to $SL_n(q)$.

As another application of our results, we classify for the first time complex representations of 
$SL_{n}(q)$ whose reductions modulo $\ell$ are 
irreducible, relying on a similar result for $GL_{n}(q)$ from 
\cite{JM} (as above, we are assuming  $\ell{\not{|}}q$). 

For a finite group $G$, let $\operatorname{IBr}(G)$ (resp. $\operatorname{Irr}(G)$) denote  
the set of isomorphism classes of irreducible $\FF G$-modules (resp. $\CC G$-modules) or the 
set of irreducible $\ell$-Brauer characters (resp. complex characters), depending on the context. 
As a final application, we exhibit  an explicit subset of $\Irr(SL_n(q))$ of size  
$|\IBr(SL_{n}(q))|$ and a partial order on $\IBr(SL_{n}(q))$ such that the corresponding  
decomposition submatrix is lower unitriangular. 
Combining this with a recent result of Bonnaf\'e \cite{Bo}, we get a parametrization of irreducible $\FF SL_n(q)$-modules. 
This also yields an ordinary basic set 
for each $\ell$-block of $SL_{n}(q)$. 

We outline an analogy with the  alternating and symmetric groups $A_n\lhd S_n$. 
There are two 
cases to consider when one studies irreducible $\FF A_n$-modules via restrictions of 
irreducible $\FF S_n$-modules. 
The case $\ell=2$, where Clifford theory is of little help, was treated in \cite{B} using 
 double covers of $S_n$. In the case $\ell\neq2$, Clifford theory 
tells us that an irreducible $\FF S_n$-module $V$ splits into two components or remains 
irreducible on restriction to $A_n$, depending on whether $V$ is isomorphic to 
$V\otimes \operatorname{sgn}$ or not, where $\operatorname{sgn}$ is the sign representation of 
$S_n$. The first explicit description of $V\otimes\operatorname{sgn}$  in terms of partitions has been found in \cite{BrIII,FK}.  As for the problem of irreducible reductions modulo $\ell$, we refer the reader to \cite{Fayers} and references therein (for alternating groups a complete classification is still only conjectural). 

\subsection{Statement of the main results}\label{SStat}
For a partition $\la=(\la_1,\la_2,\dots)$, denote  $|\la|:=\la_1+\la_2+\dots$, and 
write $\la'$ for the transposed partition. Set $\Delta(\la):=\gcd(\la_1,\la_2,\dots)$. 
For a multipartition $\underline{\la}=(\la^{(1)},\la^{(2)},\dots,\la^{(a)})$ (which 
means that each $\la^{(i)}$ is a partition) we write $\underline{\la}'$ for the transposed 
multipartition $((\la^{(1)})',\dots,(\la^{(a)})')$, and set 
$$\Delta(\underline{\la}):=\gcd(\Delta(\la^{(1)}),\dots,\Delta(\la^{(a)})).$$

For $\si\in\bar{\FF}_q$ we denote by $[\si]$ the set of all roots of the minimal polynomial of 
$\si$ over $\FF_q$; in particular, $\#[\si] = \deg(\si)$. 
We say that $\si_1$ and $\si_2$ are {\em conjugate} if $[\si_1]=[\si_2]$. 
The {\em order} of $\si$, denoted $|\si|$, 
is its multiplicative order, and $\si$ is an $\ell$- (resp. $\ell'$-) element if $|\si|$ is a 
power of $\ell$ (resp. prime to $\ell$). If $\ell=0$, all elements are $\ell'$-elements.

We state a classification of irreducible $\FF GL_n(q)$-modules, 
referring the reader to $\S$\ref{SRTGL} for details and references. 
An {\em $n$-admissible tuple} is a tuple 
\begin{equation}\label{EAT}
\big(([\si_1],\la^{(1)}),\dots,([\si_a],\la^{(a)})\big)
\end{equation}
 of pairs, where $\si_1,\dots,\si_a\in\bar{\FF}_q$ are $\ell'$-elements, and 
$\la^{(1)},\dots,\la^{(a)}$ are partitions, such that $[\si_i]\neq[\si_j]$ for all 
$i\neq j$ and $\sum_{i=1}^a \deg(\si_i) \cdot |\la^{(i)}|=n$. 
An equivalence class of the $n$-admissible tuple (\ref{EAT}) up to a permutation of the 
pairs $([\si_1],\la^{(1)}),\dots,([\si_a],\la^{(a)})$ is called an {\em $n$-admissible symbol} and denoted 
\begin{equation}\label{EAS}
\sfr 
= \big[([\si_1],\la^{(1)}),\dots,([\si_a],\la^{(a)})\big]
\end{equation}
The set ${\mathfrak L}$ of $n$-admissible symbols is the labeling set 
for irreducible $\FF GL_n(q)$-modules.
The module corresponding to the symbol (\ref{EAS}) is 
$$L({\mathfrak s}) = 
  L(\si_{1},\la^{(1)}) \circ \ldots \circ L(\si_{a},\la^{(a)}),$$ 
and 
$
\{L({\mathfrak s})\mid {\mathfrak s}\in{\mathfrak L}\}
$
is a complete set of representatives of irreducible $\FF GL_n(q)$-modules. 

The subgroup  $\OLA(\FF_q^\times)$ (consisting of the $\ell'$-elements in the multiplicative 
group $\FF_q^\times$) acts on the set ${\mathfrak L}$ of 
$n$-admissible symbols via
$$
\tau\cdot\big[([\si_1],\la^{(1)}),\dots,([\si_a],\la^{(a)})\big]=
  \big[([\tau\si_1],\la^{(1)}),\dots,([\tau\si_a],\la^{(a)})\big]
$$
for $\tau\in \OLA(\FF_q^\times)$.  
The order of the stabilizer group in  $\OLA(\FF_q^\times)$ of a symbol 
${\mathfrak s}\in {\mathfrak L}$ is called the {\em $\ell'$-branching number}
 of ${\mathfrak s}$ and is denoted 
$\kappa_{\ell'}({\mathfrak s})$. Next, if ${\mathfrak s}$ is of the form (\ref{EAS}), 
then the {\em $\ell$-branching number}  
$\kappa_{\ell}({\mathfrak s})$ is the $\ell$-part of 
$\gcd\big(n,q-1,\Delta(\underline{\la}')\big)=\gcd\big(q-1,\Delta(\underline{\la}')\big)$, 
where $\underline{\la}=(\la^{(1)},\la^{(2)},\dots,\la^{(a)})$.

\begin{Theorem}\label{TMain}
Let $V = L({\mathfrak s})$ be 
the irreducible $\FF GL_n(q)$-module corresponding to 
${\mathfrak s}\in {\mathfrak L}$. Then $V\dar_{SL_n(q)}$ is a sum of 
 $\kappa_{\ell'}({\mathfrak s})\cdot\kappa_{\ell}({\mathfrak s})$ irreducible summands.
\end{Theorem}

By Lemma~\ref{tensor}, for ${\mathfrak s}$ and ${\mathfrak s}'$ in the same 
$\OLA(\FQT )$-orbit on ${\mathfrak L}$, the restrictions 
$L({\mathfrak s})\dar_{SL_n(q)}$ and $L({\mathfrak s}')\dar_{SL_n(q)}$ are isomorphic. 
Moreover, by Theorem~\ref{TMain}, 
we have a decomposition into non-isomorphic irreducible components: 
$$L({\mathfrak s})\dar_{SL_n(q)}=
  \bigoplus_{j=1}^{\kappa_{\ell'}({\mathfrak s})\kappa_{\ell}({\mathfrak s})} L({\mathfrak s})_j$$
From Lemma~\ref{LGain} we now get:

\begin{Corollary}\label{CLabel}
The set 
$
\{L({\mathfrak s})_j\},
$
where ${\mathfrak s}$ runs through the $\OLA(\FQT )$-orbit representatives on 
${\mathfrak L}$ and $j$ runs through the integers between $1$ and 
$\kappa_{\ell'}({\mathfrak s})\kappa_{\ell}({\mathfrak s})$, is a complete set of representatives of the 
irreducible $\FF SL_n(q)$-modules. 
\end{Corollary}

Note that our labeling of the summands in the corollary is not canonical. This problem will be resolved with the use of Theorem~\ref{main3} below.

To state the result on 
reductions modulo ${\mathbf\ell}$, we need more notation. 
Given two partitions 
$\la = (\la_1,\la_2,\dots)$ and $\mu = (\mu_1,\mu_2,\dots)$ and $k \in\ZZ_{> 0}$, define new partitions  
$\la+\mu:=(\la_1 + \mu_1,\la_2+\mu_2,\dots)$ and by $k\la:=(k\la_1,k\la_2,\dots)$.

Assume until the end of the section that $\ell>0$ and as usual $\ell{\not{|}}q$. All  representations in characteristic $0$ will be over $\CC$, and the corresponding index will be used to distinguish from the modular representations. Thus the labeling set for 
irreducible $\CC GL_n(q)$-modules is the set ${\mathfrak L}_{\CC}$ of $n$-admissible symbols 
of form (\ref{EAS}), where the $\si_{i}$ need not be $\ell'$-elements, but 
$[(\si_{i})_{\ell'}] \neq [(\si_{j})_{\ell'}]$ whenever $i \neq j$. 
The irreducible $\CC GL_n(q)$-module corresponding 
to the $n$-admissible symbol (\ref{EAS}) is  
$$\SC(\sfr) = 
  \SC(\si_{1},\la^{(1)}) \circ \ldots \circ \SC(\si_{a},\la^{(a)}).$$

Given $\sfr\in{\mathfrak L}_{\CC}$ as in (\ref{EAS}), we define $\sfrs \in{\mathfrak L}$ as follows. 
For each $i$, let $s_{i}$ 
(resp. $u_{i}$) be the $\ell'$-part (resp. the $\ell$-part) of $\si_{i}$, and let 
$d_{i} := \deg(s_{i})$, $k_{i} := \deg(\si_{i})/\deg(s_{i})$. Permuting if necessary, we may assumethat $s_{1}, \ldots ,s_{m}$ form a complete set of
representatives of the conjugacy 
classes $[s_{i}]$, $1 \leq i \leq a$ (one for each). 
For $1\leq j\leq m$, set   
$$\{i_{1}, \ldots ,i_{b_j}\}=\{i\mid [s_{i}]= [s_{j}]\},$$ and   
$\delta^{(j)} := \left( k_{i_{1}}(\la^{(i_{1})})' + k_{i_{2}}(\la^{(i_{2})})' + \ldots 
  + k_{i_{b_j}}(\la^{(i_{b_j})})' \right)'
$. Now 
$$\sfrs 
:= \big[([s_{1}],\delta^{(1)}), \ldots ,([s_{m}],\delta^{(m)})\big].
$$
Finally, we say that the symbol $\sfr$ is {\em critical}, if $\ell = 2$, $q \equiv 3 (\mod 4)$,
$k_{i} > 1$ for all $i$, and there exists $j$ such that $k_{j} = 2$,
$d_{j}$ is odd, and $|u_{j}| \geq 8$.

\begin{Theorem}\label{main2}
Let $V_\CC = \SC(\sfr)$ be an irreducible $\CC GL_{n}(q)$-module and 
$W_\CC$ be an irreducible constituent of $(V_\CC)\dar_{SL_{n}(q)}$. Then reduction
modulo $\ell$ of $W_\CC$ is irreducible 
if and only if the following conditions hold:

{\rm (1)} reduction modulo $\ell$ of $V_\CC$ is irreducible 
(so $V_\CC$ is known \cite[(4.33)]{JM}); 

{\rm (2)} $\kappa_{\ell'}(\sfr) = \kappa_{\ell'}(\sfrs)$; and 
 
{\rm (3)} $\sfr$ is not critical.\\   
In this case the reduction is isomorphic to 
an irreducible constituent of $L(\sfrs) \dar_{SL_{n}(q)}$.     
\end{Theorem}

We refer the reader to Example~\ref{E1} for an illustration. 

Next, 
we show that a certain submatrix of the $\ell$-decomposition matrix  of $SL_{n}(q)$ is unitriangular. To be more precise, in $\S$\ref{SDec} we construct an explicit injective map $\Theta~:~ \IBr(SL_{n}(q)) \to \Irr(SL_{n}(q))$ and a 
partial order $\tri$ on $\IBr(SL_{n}(q))$ such that: 

\begin{Theorem}\label{main3}
For any 
$\varphi \in \IBr(SL_{n}(q))$, 
$$\Theta(\varphi) (\mod \ell) = \sum_{\psi \in \IBr(SL_{n}(q))}d_{\varphi,\psi}\psi$$
with $d_{\varphi,\varphi} = 1$, $d_{\varphi,\psi} \in \ZZ_{\geq 0}$, and
$d_{\varphi,\psi} = 0$ unless $\psi \tri \varphi$. 
\end{Theorem}

Note that Bonnaf\'e \cite{Bo} has obtained a {\em canonical} parametrization of the summands $L_\CC(\sfr)_j$ of $L_\CC(\sfr)\dar_{SL_n(q)}$, which is independent of $\sfr$ and dependent only on a fixed choice of a Gelfand-Graev representation of $SL_n(q)$, thus obtaining a parametrization of irreducible $\CC SL_n(q)$-modules. Combining this with Theorem~\ref{main3} yields a parametrization of irreducible $\FF SL_n(q)$-modules compatible with branching from $GL_n(q)$, see Proposition~\ref{cf-sl} for more details.

\noindent{\bf Acknowledgements.} We are grateful to C\'edric Bonnaf\'e for useful explanations. 



\section{Preliminaries}\label{SGen}
\subsection{Notation}
Throughout the paper all groups are assumed to be finite. We use the following notation in addition to the one introduced in $\S$\ref{SIntro}:
\begin{enumerate}
\item[$\la\vdash k$] means that $\la$ is a partition of $k$;
\item[$\tri$] denotes the dominance order on partitions;
\item[$m|\la$] means that $m$ divides all parts of the partition $\la$ (for $m\in\ZZ,\la\vdash k$);
\item[$\underline{\la}\vdash\underline{k}$] means that $\underline{\la}=(\la^{(1)},\dots,\la^{(a)})$ a multipartition of $\underline{k}=(k_1,\dots,k_a)$, i.e. each $\la^{(i)}$ is a partition of $k_i$;
\item[$g^G$] the conjugacy class of an element $g$ in a group $G$;
\item[$\overline{V}$]$=V(\mod\ell)$ reduction modulo $\ell$ of a $\CC G$-module $V$;
\item[$\eps$] a fixed generator of the cyclic group $\FF_q^\times$;
\item[$GL_n$] denotes $GL_n(q)$ and $SL_n$  denotes $SL_n(q)$.
\end{enumerate}
If $X \leq Y$ are groups, $U$ an
$\FF X$-module, and $V$ an $\FF Y$-module, then $V\dar_{X}$ is $V$ restricted to $X$, and $U\uar^{Y}$ is $U$ induced to $Y$.
Denote by $\kappa_X^Y(V)$ the number of irreducible components of $V\dar_X$ (sometimes called a {\em branching number}). We also use the notation $\kappa_X(V), \kappa_X, \kappa(V)$ or just $\kappa$ when the missing ingredients are clear from the context.

Finally, let $R_n=R_n(q)$ be the subgroup such that $SL_n\leq R_n\leq GL_n$ and 
$$R_n/SL_n=\OLA(GL_n/SL_n).$$



\subsection{Conjugacy classes in $\mathbf{GL_n}$}\label{SSGL}
If $\si\in \bar \FF_q^\times$ is an element of 
degree $d$ over $\FF_q$ then $\{1,\si,\dots,\si^{d-1}\}$ is an $\FF_q$-basis of $\FF_q(\si)= \FF_{q^d}$. The left multiplication by an element of $\FF_{q^d}$ is $\FF_q$-linear. In this way we get a group embedding $\iota^\si:\FF_{q^d}^\times\to GL_d(q)$. 
This obviously generalizes to an embedding 
$$
\iota^\si_k:GL_k(q^d)\to GL_{kd}(q)\qquad (k\in\ZZ_{>0}).
$$

Semisimple conjugacy classes of $GL_n$ are represented (uniquely up to  block permutation) by the block-diagonal matrices of the form 
\begin{equation}\label{ESSE}
s=\operatorname{diag}((\iota^{\si_1}(\si_1))^{k_1},\dots,(\iota^{\si_a}(\si_a))^{k_a}),
\end{equation}
where elements $\si_1,\dots,\si_a\in\bar \FF_q^\times$ of 
degrees $d_1,\dots,d_a$ respectively are not conjugate to each other, and  $n=k_1d_1+\dots+k_ad_a$. Let $\mathcal{S}$ be a set of representatives of the semisimple classes of $GL_n$ which are of the form (\ref{ESSE}) and $\mathcal{S}_{\ell'}$ be the subset of the  $\ell$-regular elements in $\mathcal{S}$.
For $s\in\mathcal{S}$ set
$$
\underline{k}(s):=(k_1,\dots,k_a).
$$
The centralizer $C_{GL_n(q)}(s)$ equals the (not necessarily split) Levi subgroup 
$$
GL_{k_1}(q^{d_1})\times\dots\times GL_{k_a}(q^{d_a}),
$$
embedded into $GL_{k_1d_1}(q)\times\dots\times GL_{k_ad_a}(q)<GL_n(q)$ via 
$\iota^{\si_1}_{k_1}\times\dots\times \iota^{\si_a}_{k_a}$. 

The unipotent conjugacy classes of $GL_n(q)$ are represented (uniquely) by the Jordan matrices
$$
J(\la):=\diag(J(\la_1),\dots,J(\la_m))\qquad(\la\vdash n),
$$
where $J(\la_i)$ is the Jordan block of size $\la_i$. If we write $\la$ in the form $\la=(1^{r_1},2^{r_2},\dots)$, then by \cite{Spr} there exists $N=N(\la)\in\ZZ_{\geq 0}$ such that 
\begin{equation}\label{ESpr}
|C_{GL_n(q)}(J(\la))|=q^N\prod_{i\geq 1}|GL_{r_i}(q)|.
\end{equation}

By the Jordan decomposition, the conjugacy classes of $GL_n$ are labeled by pairs $(s,\underline{\la})$, where $s\in\mathcal{S}$ and $\underline{\la}$ is a multipartition of $\underline{k}(s)$. 
A canonical representative of the corresponding conjugacy class looks like $su=us$ for the unipotent element $u\in GL_n$ chosen as follows. For each $i=1,\dots,a$, let $u_i$ be the Jordan matrix matrix $J(\la^{(i)})\in GL_{k_i}(q^{d_i})$, and let 
$$
u:=\diag(\iota^{\si_1}_{k_1}(u_1),\dots,\iota^{\si_a}_{k_a}(u_a))\in GL_{k_1d_1}(q)\times\dots\times GL_{k_ad_a}(q)<GL_n(q).
$$
Note that the element $su$ is $\ell$-regular if and only if $s$ is $\ell$-regular.

\subsection{Representation theory of $\mathbf{GL_n}$}\label{SRTGL}
If $n=kd$, then to {\em every}, not necessarily $\ell$-regular, 
element
$\sigma$ of degree $d$ over $\FF_q$ and $\la\vdash k$ we
associate the irreducible $\FF GL_n$-module denoted $L(\sigma,\la)$, see \cite[(3.5.3)]{BDK} (or \cite{J}, where it is denoted $D_\FF(\si,\la)$). 
Let $\circ$ denote 
the Harish-Chandra induction. 

For $s\in\mathcal{S}$ as in (\ref{ESSE}) and 
$
\underline{\la}=(\la^{(1)},\dots,\la^{(a)})\vdash \underline{k}(s)
$, define
$$
L(s,\underline{\la}):=L(\si_1,\la^{(1)})\circ\dots\circ  L(\si_a,\la^{(a)}).
$$

\begin{Theorem}\label{TIrrClass}
The set $\{L(s,\underline{\la})\mid s\in \mathcal{S}_{\ell'},\underline{\la}\vdash\underline{k}(s)\}$, is a complete set of representatives of isomorphism classes of the  irreducible 
$\FF GL_n$-modules.
\end{Theorem}

To connect with the notation used in the Introduction, define an $n$-admissible symbol  
$${\mathfrak s}:=
[([\si_1],\la^{(1)}),\dots,([\si_a],\la^{(a)})].
$$ 
Then $L(s,\underline{\la})=L({\mathfrak s})$, cf.   $\S$\ref{SStat}.
This is the James' classification of irreducible $\FF GL_n$-modules (see \cite{J} and \cite[4.4b]{BDK}). 
Even though this classification suggests that the modules $L(\si,\la)$ for $\ell$-singular $\si$ 
are redundant, it is sometimes convenient to use them. So we will not assume that $\si$ is $\ell$-regular, 
unless otherwise stated. 

If $\FF=\CC$, we write $L_\CC(\si,\la)$ instead of $L(\si,\la)$, $L_\CC(\mathfrak{s})$ instead of $L(\mathfrak{s})$, etc. in order to distinguish from the case of positive characteristic. The {\em Specht module} $S_\FF(\si,\la)$ can be constructed in any characteristic, and $S_{\FF}(\si,\la)$ is a reduction modulo $\ell$ of $L_\CC(\si,\la)=S_\CC(\si,\la)$, see \cite{J}.

\subsection{Conjugacy classes in the intermediate subgroup}
We will make use of the following general:

\begin{Lemma}\label{LConjCl}
Let $S \lhd G$ and $S \leq R \leq G$ be such that $G/S$ is cyclic. For any
$g \in G$, set $c = (G:C_{G}(g)S)$ and $d = (G:R)$. Then
$|g^G|/|g^R| = \gcd(c,d)$. Moreover, if $G/R$ is an $\ell$-group and $g$ is
an $\ell'$-element, then $g \in R$.
\end{Lemma}
\begin{proof}
Denote $C:=C_G(g)$ and $D:=C\cap R$. Then $C\cap S=D\cap S$, hence 
\begin{equation}\label{EGain}
\frac{|DS|}{|CS|}=\left(\frac{|D||S|}{|D\cap S|}\right):\left(\frac{|C||S|}{|C\cap S|}\right)=\frac{|D|}{|C|}.
\end{equation}
Let $G/S=\langle a\rangle$ be of order $m$, and $\pi:G\to\langle a\rangle$ be a surjection with kernel $S$. 
Then $\pi(R)=\langle a^d\rangle$ and $\pi(CS)=\langle a^c\rangle$, hence $\pi(CS\cap R)=\langle a^{\lcm(c,d)}\rangle$. It is easy to see that $CS\cap R=DS$, so $|\pi(DS)|=m/\lcm(c,d)$. Now we have 
$$\frac{|g^G|}{|g^R|}=\frac{|G|}{|C|}:\frac{|R|}{|D|}=\frac{|G|}{|R|}\frac{|D|}{|C|}=\frac{|G|}{|R|}\frac{|DS|}{|CS|}=\frac{|G|}{|CS|}\frac{|DS|}{|R|}=c\frac{m/\lcm(c,d)}{m/d},
$$
which is $\gcd(c,d)$, as required.
\end{proof}

To apply Lemma~\ref{LConjCl}, we need the following fact.

\begin{Lemma}\label{LDet}
Let $u\in GL_n(q)$ be a unipotent element whose conjugacy class corresponds to a partition $\la$. Then $\det$ maps $C_{GL_n(q)}(u)$ onto the subgroup $\langle\eps^{\Delta(\la)}\rangle$ of $\FF_q^\times$. 
\end{Lemma}
\begin{proof}
Write $\la=(1^{r_1},2^{r_2},\dots)$ and $C:=C_{GL_n(q)}(u)$.  Then $u$ is conjugate to the matrix 
$\bigoplus_{i:r_i>0} I_{r_i}\otimes J(i)$, where $I_r$ is the identity $r\times r$ matrix.
Hence $C$ contains the subgroup $D:=\prod_{i:r_i>0} GL_{r_i}(q)\otimes I_i$. Let $P$ be a Sylow $p$-subgroup of $C$. 
Using (\ref{ESpr}), we conclude that $C=PD$. Note that $\det(P)$ is trivial, so $\det(C)=\det(D)$. Now, for each $i$, $\det$ maps the subgroup $GL_{r_i}(q)\otimes I_i$ onto $\langle\eps^{i}\rangle$. So it maps $D$ onto $\prod_{i:r_i>0}\langle\eps^{i}\rangle=\langle\eps^{\Delta(\la)}\rangle$. 
\end{proof}


\begin{Proposition}\label{PIndex}
Let $g=su$ be a standard representative of an $\ell$-regular conjugacy class in $GL_n$ corresponding to the pair $(s,\underline{\la})$, where $s\in\mathcal{S}_{\ell'}$ 
and $\underline{\la}=(\la^{(1)},\dots,\la^{(a)})$ is a multipartition of $\underline{k}(s)$. 
Then $g\in R_n$ and 
$$
\frac{|g^{GL_n}|}{|g^{R_n}|}=\gcd\big\{(GL_n:R_n),\Delta(\underline{\la})\big\}.
$$ 
\end{Proposition}
\begin{proof}
In view of Lemma~\ref{LConjCl}, we need to know the image of $C:=C_{GL_n(q)}(g)$ under the determinant map $\det:GL_n(q)\to\FF_q^\times$. 
Now 
\begin{equation}\label{CSU}
C=C_{GL_n(q)}(su)=C_{C_{GL_n(q)}(s)}(u)=\prod_{i=1}^a C_{GL_{k_i}(q^{d_i})}(u_i),
\end{equation}
where $u_i$ are as in $\S$\ref{SSGL}. 

For positive integers $m,d$ and an element $\sigma$ of degree $d$ over $\FF_q$, we consider the composition map $\det_{d}$:
$$
{\textstyle\det_d}:GL_m(q^d)\stackrel{\iota_m^\sigma}{\longrightarrow}GL_{md}(q)\stackrel{\det}{\longrightarrow}\FF_q^\times.
$$
Observe that $\det_d$ is also the composition
$$
{\textstyle\det_d}:GL_m(q^d)\stackrel{\det}{\longrightarrow}\FF_{q^d}^\times\stackrel{N_{\FF_{q^d}/\FF_q}}{\longrightarrow}\FF_q^\times.
$$
(In particular, $\det_d$ does not depend on $\si$). Fix a generator $\eps_d$ of the group $\FF_{q^d}^\times$ such that $\eps=N_{\FF_{q^{d_i}}/\FF_q}(\eps_d)$. 
By Lemma~\ref{LDet}, $\det$ maps $C_{GL_{k_i}(q^{d_i})}(u_i)$ onto $\langle\eps_{d_i}^{\Delta(\la^{(i)})}\rangle$, and  $N_{\FF_{q^{d_i}}/\FF_q}$ maps the latter onto $\langle\eps^{\Delta(\la^{(i)})}\rangle$. So the image of $C$ under $\det$ is $\langle\eps^{\Delta(\underline{\la})}\rangle$. The result now follows from Lemma~\ref{LConjCl}. 
\end{proof}

\section{Clifford Theory}
\subsection{Known results}
Let $S\lhd G$ be a normal subgroup and $V\in\IBr(G)$. Clifford theory (see e.g.
\cite[III.2]{F}, \cite[$\S$49]{CR}, \cite[$\S$11A]{CRI}) gives information on the restriction $V\dar_S$. For example, it is known that  
$V\dar_S$ is completely reducible, so we can write $$V \dar_{S} = \bigoplus^{t}_{i=1}V_{i},\qquad  (V_{1}, \ldots ,V_{t} \in \IBr(S)).$$ 
Moreover, the irreducible $S$-modules $V_1,\dots, V_t$ are  
$G$-conjugate to each other (but not necessarily pairwise non-isomorphic in general).
Let $I: = \Stab_{G}(V_{1})$ be the {\em inertia group} of $V_{1}$ in $G$. Then $I\geq S$ and there is an $\FF I$-module $\tilde V_1$ such that $V_1$ is an irreducible component of $\tilde V_1\dar_S$ and  $V = \tilde V_1 \uar ^{G}$.

If $G/S$ is an $\ell$-group and an irreducible $\FF S$-module appears as a component of the restrictions $V\dar_S$ and $V'\dar_S$ of the irreducible $\FF G$-modules $V$ and $V'$, then $V\cong V'$.
If $G/S$ is cyclic, $V_{1}, \ldots ,V_{t}$ are pairwise non-isomorphic, $\tilde V_1=V_1$ as 
vector spaces, and 
$\kappa^{G}_{S}(V)=t= (G:I)$. 

\begin{Lemma}\label{Llink2}
Let $S\leq H$ be normal subgroups of a group $G$, $V\in\IBr(G)$, and $U$ be an irreducible component of $V\dar_H$. Then $\kappa^{G}_{S}(V)=\kappa^{G}_{H}(V)\cdot \kappa^{H}_{S}(U)$.
\end{Lemma}
\begin{proof}
Write $V \dar_{S} = \oplus^{t}_{i=1}V_{i}$, 
where $V_{1}, \ldots ,V_{t} \in \IBr(S)$ are $G$-conjugate, and 
$\kappa^{G}_{S}(V) = t$. 
Similarly, write $V \dar_{H} = \oplus^{m}_{j=1}U_{j}$, 
where $U=U_{1}, \ldots ,U_{m} \in \IBr(H)$ are $G$-conjugate, and 
$\kappa^{G}_{H}(V) = m$. 
We may assume 
that $V_{1}$ is a 
constituent of $U_1\dar_{S}$. Then 
$U_{1} \dar_{S} = \oplus^{a}_{j=1}U_{1,j}$, 
where $V_{1} \simeq U_{1,1}, \ldots ,U_{1,a} \in \IBr(S)$ are $H$-conjugate, and 
$a = (H:H_{1})$ for $H_{1} = \Stab_{H}(V_{1})$. For any $i$ there is $g \in G$ such that 
$V_{i} \simeq {}^gV_{1}$. Hence
$$\Stab_{H}(V_{i}) = H \cap g^{-1} \Stab_{G}(V_{1})  g
=g^{-1}(H \cap  \Stab_{G}(V_{1}))g =g^{-1}H_{1}g$$
since $H\lhd G$. Thus $\kappa^{H}_{S}(U_{j}) = a = \kappa^{H}_{S}(U)$ for all $j$, 
and so $t = ma$.
\end{proof}

\subsection{Number of constituents for cyclic extensions}
We need more results in spirit of Clifford theory, especially for the case of cyclic extensions, i.e. when $G/S$ is cyclic. Given a $G/S$-module 
we consider it as a $G$-module via inflation without further comment.

\begin{Lemma}\label{linear}
 Let $S\lhd G$ with $G/S$ a cyclic $\ell'$-group, and  
$V \in \IBr(G)$.
\begin{enumerate}
\item[{\rm (i)}] We have $\kappa_S(V) = \sharp\{L \in \IBr(G/S) \mid V \simeq V \otimes L\}$.

\item[{\rm (ii)}] If $L \in \IBr(G/S)$ and 
$V \simeq V \otimes L$ then $\kappa_{S}(V)\geq (G:\Ker(L))$.
\end{enumerate}
\end{Lemma}

\begin{proof}
(i) Let $\eta:=\sharp\{L \in \IBr(G/S) \mid V \simeq V \otimes L\}$. If $L \in \IBr(G/I)$, we have 
$$V \otimes L = (\tilde V_1 \uar ^{G}_I) \otimes L \simeq (\tilde V_1 \otimes (L\dar_{I})) 
  \uar ^{G} \simeq \tilde V_1 \uar ^{G} =  V,$$
and so $\eta \geq \kappa_S(V)$. 
Conversely, assume that $V \simeq V \otimes L$ for some 
$L \in \IBr(G/S)$. To see that $\eta \leq \kappa_S(V)$, it suffices to show that $I$ acts trivially on $L$. 
Again we have 
$\tilde V_1 \uar^G \simeq V \simeq V \otimes L \simeq (\tilde V_1 \otimes (L\dar_{I})) \uar^{G}.$ 
As a subspace $\tilde V_1$ equals the  
$V_{1}$-isotypic component of the $S$-module $V$. 
It follows that
$\tilde V_1 \simeq \tilde V_1 \otimes (L\dar_{I})$, so  
$\tilde V_1^{*} \otimes \tilde V_1 \simeq \tilde V_1^{*} \otimes \tilde V_1 \otimes (L\dar_{I})$. As $S$ is irreducible on
$\tilde V_1$ and trivial on $L$, taking $S$-fixed points, we see that $1_{I} \simeq L \dar_{I}$, as required.

(ii) Arguing as in (i), we see that $L\dar_{I}$ is trivial. So $\kappa_S(V) = (G:I) \geq (G:\Ker(L))$. 
\end{proof}

Lemma \ref{linear} 
provides information on $\kappa(V)$ 
for 
cyclic $\ell'$-extensions. 
For an {\em arbitrary} 
cyclic extension 
the problem breaks into two parts; 
the first (easier) part for a cyclic $\ell'$-extension, and the second 
part for a cyclic $\ell$-extension:

\begin{Lemma}\label{rami}
{ Let $r$ be a prime, $S \lhd G$ with $G/S$ cyclic, 
$V\in\IBr(G)$, and $S \leq A,B \leq G$ be such that $A/S = O_{r}(G/S)$ and  
$B/S = O_{r'}(G/S)$. Also, let $U$ (resp. $W$) be an irreducible constituent of $V\dar_{A}$ 
(resp. $V\dar_{B}$). Then
\begin{enumerate}
\item $\kappa^{G}_{S}(V) = \kappa^{G}_{A}(V) \cdot \kappa^{G}_{B}(V)$;

\medskip
\item $\kappa^{G}_{A}(V) = \kappa^{B}_{S}(W)$, $\kappa^{G}_{B}(V) = \kappa^{A}_{S}(U)$.
\end{enumerate}}
\end{Lemma}
 
\begin{proof}
By Lemma~\ref{Llink2}, 
$$t:=\kappa^{G}_{S}(V) = \kappa^{G}_{A}(V)\kappa^{A}_{S}(U)= \kappa^{G}_{B}(V)\kappa^{B}_{S}(W).$$ Write $V \dar_{S} = \oplus^{t}_{i=1}V_{i}$ for $V_{1}, \ldots ,V_{t} \in \IBr(S)$, and let $I = \Stab_{G}(V_{1})$, $A_1=\Stab_{A}(V_{1})$, $B_1=\Stab_{B}(V_{1})$. 
Then $G/S = A/S \times B/S$ and $I/S \geq A_{1}/S \times B_{1}/S$. Conversely, 
assume $x \in I$ and write $xS = yS \cdot zS$, where $yS$ is the $r$-part
of $xS$ and $zS$ is the $r'$-part of $xS$. Then $yS$ and $zS$ are powers of $xS$ and 
so $y$ and $z$ stabilize $V_{1}$, hence $y \in A_{1}$ and $z \in B_{1}$. 
Thus $I/S = A_{1}/S \times B_{1}/S$. So 
$$t = (G:I) = (A:A_{1}) \cdot (B:B_{1}) = \kappa^{A}_{S}(U)\kappa^{B}_{S}(W).$$ 
Hence $\kappa^{B}_{S}(W) = \kappa^{G}_{A}(V)$, 
$\kappa^{A}_{S}(U) = \kappa^{G}_{B}(V)$, and
$\kappa^{G}_{S}(V)  = \kappa^{G}_{A}(V) \cdot \kappa^{G}_{B}(V)$.
\end{proof}

Next we link branching numbers $\kappa^{G}_{S}$  
to those for intermediate extensions $S \lhd H$ and $H \lhd G$.

\begin{Lemma}\label{link1}
Let $S\lhd G$ with $G/S$ a cyclic $\ell$-group, $S\leq H\leq G$, $V \in \IBr(G)$. 
Then 
 $\kappa^{G}_{H}(V) = \min\{ \kappa^{G}_{S}(V), |G/H|\}$.
\end{Lemma}
 
\begin{proof}
Let $U \in \IBr(H)$ an irreducible constituent of $V\dar_{H}$ and $W \in \IBr(S)$ be an irreducible constituent of $U\dar_{S}$.
Let $I := \Stab_{G}(W)$. 
Then $V = \tilde W\uar^{G}$ for $\tilde W \in \IBr(I)$ with $\tilde W\dar_{S} = W$.
As $G/S$ is a cyclic $\ell$-group, either $I \geq H$ or $H > I$.
If $I \geq H$ then $U' := \tilde W\dar_{H}$ is also irreducible and 
$I$-invariant. Conversely, if $g \in G$ stabilizes $U'$ then $g$ stabilizes
$U'\dar_{S} = W$, so $g \in I$. Thus $\Stab_{G}(U') = I$, whence 
$$\kappa^{G}_{H}(V) = (G:I) = \kappa^{G}_{S}(V) \leq |G/H|,$$ 
and we are done in the case $I \geq H$.

If $H > I$ then $U' := \tilde W\uar^{H}$ is an irreducible component of $V\dar_{H}$. Hence 
$\kappa^{G}_{H}(V) = \dim V/\dim U'  = |G/H| < (G:I),$ 
and we are done again.
\end{proof}

\begin{Corollary}\label{link2}
Let $S\lhd G$ with $G/S$ cyclic, $S\leq H\leq B\leq  G$ such that $B/H = \OLB(G/H)$, 
and $V \in \IBr(G)$. 
Then
$$\kappa^{G}_{H}(V) = \kappa^{G}_{B}(V) \cdot \min\{ (\kappa^{G}_{S}(V))_{\ell}, |G/H|_{\ell}\}.$$
\end{Corollary}

\begin{proof}  
Consider $A \geq H$ with $A/H = \OLA(G/H)$, and $R,T \geq S$ with 
$R/S = \OLA(G/S)$, $T/S = \OLB(G/S)$. By Lemma \ref{rami}, 
$\kappa^{G}_{H}(V) = \kappa^{G}_{A}(V) \cdot \kappa^{G}_{B}(V)$ so we just need to compute $ \kappa^{G}_{A}(V)$. 
As 
$$RH/H \simeq R/(R\cap H) \simeq (R/S)/((R \cap H)/S)$$ 
is 
an $\ell'$-group, $RH/H \leq \OLA(G/H) = A/H$, whence $R \leq A$. By  Lemma~\ref{link1}, $\kappa^{G}_{A}(V) = \min\{ \kappa^{G}_{R}(V), |G/A|\}$. 
Clearly $|G/A| = |G/H|_{\ell}$. By Lemma~\ref{rami}, 
$\kappa^{G}_{S}(V) = \kappa^{G}_{R}(V) \cdot \kappa^{G}_{T}(V)$. As the first factor
is an $\ell$-power and the second factor is coprime to $\ell$, we have 
$\kappa^{G}_{R}(V) = (\kappa^{G}_{S}(V))_{\ell}$, and  we are done.  
\end{proof}  

Corollary~\ref{link2} shows that if we know 
$\kappa^G_S(V)$, then, to 
determine $\kappa^G_H(V)$ for $H \geq S$, it suffices to know $\kappa^G_B(V)$, 
which can be found using Lemma~\ref{linear}. 


\begin{Lemma}\label{cross}
 Let $G = G_1 \times G_2$, $V \in \IBr(G)$, $H_{i} \lhd G_i$ and 
$G_i/H_{i} = \langle \bar{x_i} \rangle \simeq C_{n}$ for $i =
 1,2$. Let the subgroup $K\leq G$ 
be such that $G \geq K \geq H_{1} \times H_{2}$ and 
$K/(H_{1} \times H_{2}) = \langle \bar{x_1}\bar{x_2} \rangle$. Then
$\kappa_{K}(V)=\gcd(s_1,s_2)$,
where $s_1:=\kappa_{H_{1} \times G_2}(V)$, $s_2:=\kappa_{G_1 \times H_{2}}(V)$.
\end{Lemma}

\begin{proof}
Throughout the proof we assume  $i\in\{1,2\}$. Write $V = U^1 \otimes U^2$ for $U^i \in \IBr(G_i)$. By assumption, 
$U^i\dar_{H_{i}} = U^i_{1} \oplus \ldots \oplus U^i_{s_i}$ for  $U^i_{1}, \ldots ,U^i_{s_i} \in \IBr(H_{i})$ 
and $s_i = (G_i:I_{i})$ where $I_{i} := \Stab_{G_i}(U^i_{1})$. 
Clearly, $J := \Stab_{K}(U^1_{1} \otimes U^2_{1})$ contains $H_{1} \times H_{2}$. Choose
$x_i \in G_i$ such that $\bar{x_i} = x_iH_{i}$. Then 
$K = \langle H_{1},H_{2},x_1x_2 \rangle$ and $G/K \simeq C_{n}$. 
Assume that $(x_1x_2)^{m} \in J$ for some integer $m$.
Then 
$$U^1_{1} \otimes U^2_{1} \simeq (x_1^{m}x_2^{m})(U^1_{1} \otimes U^2_{1}) \simeq 
(x_1^{m}U^1_{1}) \otimes (x_2^{m}U^2_{1})$$
as $H_{1} \times H_{2}$-modules. By the irreducibility of the $U^i_{1}$, we see that 
$x_i^{m}U^i_{1} \simeq U^i_{1}$, i.e. $x_i^{m} \in I_{i}$.
But $G_i/I_{i}$ is cyclic of order $s_i$, hence $s_i|m$. Setting 
$r = {\mathrm {lcm}}(s_1,s_2)$, we see that $r|m$. Thus $(x_1x_2)^{m} \in J$ only if $r|m$. The converse 
is obviously true. Therefore, $J/(H_{1} \times H_{2}) = \langle (\bar{x}\bar{y})^{r} \rangle$ and 
$(K:J) = r$. 

Now choose $T$ to be an irreducible constituent of $V\dar_{K}$ that lies over 
$U^1_{1} \otimes U^2_{1}$. As 
$(K:J) = r$, we have $\kappa_{H_{1} \times H_{2}}(T)=r$. Hence
$$\dim(T) = r\dim(U^1_{1})\dim(U^2_{1}),$$ 
so 
$\dim(V)/\dim(T) = s_1s_2/r = \gcd(s_1,s_2)$, and $\kappa_{K}(V)=\gcd(s_1,s_2)$.  
\end{proof}

\begin{Lemma}
\label{2mod}
Let $S\lhd G$ with $G/S$ cyclic, $S\leq A\leq G$ with $A/S=O_\ell(G/S)$, 
and $U \in \IBr(S)$. Assume that $U$ is an $S$-composition factor of $V_{i}\dar_{S}$ for some
$V_{i} \in \IBr(G)$, $i = 1,2$. Then 
$V_{2} = V_{1} \otimes L$ for some $L \in \IBr(G/A)$.
\end{Lemma}

\begin{proof}
By assumption, there are $W_i\in\IBr(A)$ such that $U$ is a direct summand of $W_i\dar_A$ and $W_i$ is a direct summand of $V_i\dar_A$. Since $A/S$ is an $\ell$-group, $W_1\cong W_2\cong W$. 
Let $I := \Stab_{G}(W)$. By Clifford theory, we have $V_{i} = M_{i}\uar^{G}$ for some 
$M_{i} \in \IBr(I)$ and $W=M_{i}\dar_{A}$  ($i = 1,2$). It follows that 
 $M_{2} \simeq M_{1} \otimes N$ for some 
$N \in \IBr(I/A)$. But $G/A$ is cyclic, hence $N$ is also $G/A$-invariant and it extends to 
an irreducible $G/A$-module $L$. Thus 
$$V_{2} \simeq (M_{1} \otimes N)\uar^{G} = (M_{1} \otimes (L\dar_{I}))\uar^{G} \simeq 
 (M_{1} \uar^{G}) \otimes L = V_{1} \otimes L,$$
 as required.  
\end{proof}

We now address the question of describing irreducible $\FF S$-modules in terms of irreducible $\FF G$-modules and branching from $G$ to $S$. 

\begin{Lemma}\label{LGain}
Let $S\lhd G$ with $G/S$ cyclic, and $S \leq A,B \leq G$ with $A/S = \OLB(G/S)$, $B/S = \OLA(G/S)$. Introduce the equivalence relation $\sim$ on $\IBr(G)$ by setting $V\sim V'$ if and only if $V\simeq V'\otimes L$ for some  $L\in\IBr(G/A)$, and let $V_1,\dots,V_m$ be a complete set of representatives of the equivalence classes. Decompose 
$$
V_i\dar_S=\bigoplus_{j=1}^{k_i} V_{ij}\qquad(1\leq i\leq m)
$$
where $V_{ij}$ are irreducible $\FF S$-modules. Then $$\{V_{ij}\mid 1\leq i\leq m,\ 1\leq j\leq k_i\}$$ is a complete set of representatives of the isomorphism classes of the irreducible $\FF S$-modules.
\end{Lemma}
\begin{proof}
We only have to check that $V_{ij}\not\simeq V_{kl}$ if $i\neq k$, which follows immediately from Lemma~\ref{2mod}.
\end{proof}

\subsection{Clifford theory and blocks}

\begin{Lemma}
\label{1bl}
 Let $S\lhd G$ with $G/S$ an $\ell$-group, 
$V$ an $\FF G$-module, and 
\begin{enumerate}
\item[{\rm (i)}] $V\dar_{S}$ has a filtration $0 = V^{0} < V^{1} < \ldots$, where
each $V_{i} := V^{i}/V^{i-1}$ is a $G$-conjugate of 
$V_{1} = V^{1}$;


\item[{\rm (ii)}] All composition factors of $V_{1}$ belong to the same $\FF S$-block.
\end{enumerate}
Then all composition factors of $V$ belong to the same $\FF G$-block.
\end{Lemma}

\begin{proof}
Pick an $S$-composition factor $U_{1}$ of $V_{1}$ and $W_{1} \in \IBr(G)$ 
such that $U_{1}$ is a submodule of $W_{1}\dar_{S}$. Let $B$ be the 
$G$-block containing $W_{1}$, and $b$ be the $S$-block containing $U_{1}$. Then $B$ covers 
$b$, as well as any $G$-conjugate of $b$, cf. \cite[Lemma IV.4.10]{F}. By the assumption (ii), 
every composition factor of $V\dar_S$ belongs to some $G$-conjugate of $b$. 

Now consider any composition factor $M$ of $V$ and any composition factor $N$ of 
$M\dar_{S}$. Then we may assume that $N$ is a composition factor of some $V_{i}$. But 
$V_{i} = \ ^{z}V_{1}$ for some $z \in G$. Hence $N$ belongs to the $S$-block $b^{z}$. 
Now $b^{z}$ is
also covered by $B$. In fact, since $G/S$ is an $\ell$-group, $B$ is the unique $G$-block
that covers $b^{z}$, cf. \cite[Lemma V.3.5]{F}. Hence $M$ belongs to $B$.
\end{proof}

\begin{Lemma}\label{ext}
 Let $S\lhd G$ with $G/S$ cyclic, and let 
$B$ be an $\ell$-block of $S$. Suppose that every ordinary character $\rho$ belonging to 
$B$ is $G$-invariant. Then every Brauer character $\varphi$ in $B$ is extendible to $G$.
\end{Lemma}

\begin{proof}
By \cite[Theorem 61.6]{Do}, $\varphi = \sum_{\rho \in B \cap \Irr(S)}a_{\rho}\bar{\rho}$ 
for some $a_{\rho} \in \ZZ$, and $\bar{\rho}$ denotes the restriction of 
irreducible characters $\rho$ to 
$\ell'$-classes of $S$. It follows that $\varphi$ is $G$-invariant. Since $G/S$ is cyclic,
$\varphi$ is extendible to $G$.
\end{proof}

\begin{Lemma}\label{red1a}
Let $S\lhd G$ with $G/S$ a cyclic $\ell$-group, 
and 
$B$ be an $\ell$-block of $S$. Suppose there is $b \in\ZZ_{\geq 0}$
such that $(G:\Stab_{G}(\rho)) \leq \ell^{b}$ for any complex irreducible character  
$\rho \in B$. Then 
$(G:\Stab_{G}(\psi)) \leq \ell^{b}$ for any $\psi \in B \cap \IBr(S)$.
\end{Lemma}
 
\begin{proof}
We may assume that $b \leq a$, where $\ell^{a} := |G/S|$.
Let $J$ be the unique subgroup of index $\ell^{b}$ in $G$ containing $S$. By assumption,
any complex $\rho \in B$ is $J$-invariant. Now Lemma \ref{ext} applied to
$J \lhd G$ implies that every $\psi \in B \cap \IBr(S)$ is $J$-invariant. 
\end{proof}

\begin{Lemma}\label{red12}
Let $S\lhd G$ with $G/S$ a cyclic $\ell$-group,
and $B$ be an $\ell$-block of $G$. Suppose there is $b \in\ZZ_{\geq 0}$ such that $\kappa^{G}_{S}(\chi) \leq \ell^{b}$ for any complex irreducible character 
$\chi \in B$. Then 
$\kappa^{G}_{S}(\varphi) \leq \ell^{b}$ for any $\varphi \in B \cap \IBr(G)$.
\end{Lemma}
 
\begin{proof}
Let $\varphi \in B \cap \IBr(G)$, $\psi \in \IBr(S)$ be an irreducible component of 
$\varphi\dar_{S}$, and $B_{1}$ be the $\ell$-block of $S$ containing $\psi$. As $G/S$ is 
cyclic, by Lemma~\ref{red1a} it suffices to show that 
$(G:\Stab_{G}(\rho)) \leq \ell^{b}$ for any complex irreducible $\rho \in B_{1}$.
Consider any such $\rho$. As $G/S$ is an $\ell$-group, $B$ is the unique block of $G$ that 
covers $B_{1}$. So there is some irreducible $\chi \in B$ such that $\rho$ is an 
irreducible component of $\chi\dar_{S}$. As $G/S$ is cyclic, 
$(G:\Stab_{G}(\rho)) = \kappa^{G}_{S}(\chi)$, and so we are done.
\end{proof}

\subsection{${\mathbf \ell}$-elementary groups} A finite group is {\it $\ell$-elementary} if it is a direct product of
an $\ell$-group and a cyclic $\ell'$-group. It is easy to see that any subgroup and any quotient group of an $\ell$-elementary group are also $\ell$-elementary.
In this subsection we explain why many results in Clifford theory valid for cyclic extensions, are also valid  for $\ell$-elementary extensions. This material, even though 
not needed in this paper, will be useful in other situations of Clifford theory.

\begin{Lemma}\label{elem}
Any irreducible projective representation over the field $\FF$ of an $\ell$-elementary group $H$ has degree $1$.
\end{Lemma}
 
\begin{proof}
We can realize our projective representation as a linear representation $\Psi$ of a covering group $T$ of $H$. 
Writing $Z := Z(T)$, we have $T/Z \simeq H$. As $\Psi\dar_{\OLB(Z)}$
is trivial, we may assume that $Z$ is an $\ell'$-group. Now 
$\OLA(T) \geq Z$ and $\OLA(T)/Z \simeq \OLA(H)$ is cyclic. Hence 
$\OLA(T)$ is abelian. As $H$ is nilpotent, so is $T$. Hence $T = \OLA(T) \times \OLB(T)$, and $\deg(\Psi) = 1$. 
\end{proof}


\begin{Lemma}\label{rami3}
Let $S\lhd G$ with $G/S$ being $\ell$-elementary, and 
$V \in \IBr(G)$. Then $V\dar_S$ is multiplicity free.
\end{Lemma}

\begin{proof}
Let $U$ be an irreducible component of $V\dar_S$ and $I:= \Stab_{G}(U)$. 
Assume first that $G = I$. By \cite[11.20]{CRI}, 
$V = X \otimes Y$, where $X$ is a projective representation of $G$ of dimension equal to 
$\dim U$ and $Y$ is a projective representation of $G/S$. By Lemma \ref{elem}, $\dim Y=1$, so $V\dar_{S} = U$.

In the general case we have $V = \tilde U\uar^G$
for some $\tilde U \in \IBr(I)$ lying above $U$. As $U$ is $I$-invariant and  
$I/S$ is $\ell$-elementary, the above argument shows that 
$\tilde U\dar_{S} = U$. In particular, $\dim V/\dim U = (G:I)$, and the claim follows.   
\end{proof}

Now Lemmas \ref{elem} and \ref{rami3} allow us to repeat the proof of Lemma \ref{rami}.  So the conclusions of that key lemma also hold for $\ell$-elementary extensions.


\section{Main Results}


\subsection{Computing $\mathbf{\kappa_{\ell'}}$}
The irreducible $\FF GL_n$-modules which factor through $SL_n$ 
are among the $L(\tau,(n))$ with $\tau \in \OLA(\FQT )$, see \cite{J}.  
The following statement is needed in view of Lemma~\ref{linear}, 
cf. also \cite[Lemma 2.9]{GT1}:

\begin{Lemma}\label{tensor}
Let $\mathfrak{s}$ be an $n$-admissible symbol and $\tau\in \OLA(\FQT )$. Then 
$$L(\mathfrak{s}) \otimes L(\tau,(n)) \simeq L(\tau\cdot\mathfrak{s}).$$
\end{Lemma}

\begin{proof}
It suffices to prove that 
$L(\si,\la) \otimes L(\tau,(n)) \simeq L(\si\tau,\la)$, where $\deg(\si) = d$, $\la \vdash k := n/d$. We use notation and results from \cite{J}. 
First of all, the values of the Brauer character of the irreducible $\FF GL_{d}$-module 
$M_{\FF}(\si,(1))$ are given explicitly in \cite[(3.1)]{J}---they involve a certain function 
$\langle x,y\rangle_{d}~:~\FF_{q^{d}}^{\times} \to \CC^{\times}$ with the property 
$$\langle x_{1}x_{2},y\rangle_{d} = \langle x_{1},y\rangle_{d} \cdot \langle x_{2},y\rangle_{d}.$$ 
In particular, $\langle \si\tau,y\rangle_{d} = \langle \si,y\rangle_{d} \cdot \langle \tau,y\rangle_{d}$, and 
$\langle \tau,y\rangle_{d}^{q-1} = 1$. These two equalities imply that 
the Brauer character of $M_{\FF}(\si\tau,(1))$ is just the product of the Brauer characters  of 
$M_{\FF}(\si,(1))$ and $L(\tau,(n))$, whence $$M_{\FF}(\si,(1)) \otimes L(\tau,(n)) \simeq M_{\FF}(\si\tau,(1)).$$ 
This identity, together with the construction of 
$M_{\FF}(\si,(1^{k}))$ in \cite{J}, 
implies that $$M_{\FF}(\si,(1^{k})) \otimes L(\tau,(n)) \simeq M_{\FF}(\si\tau,(1^{k})).$$ 

Next, the submodule $\SF(\si,\la)$ of $M_{\FF}(\si,(1^{k}))$ is defined using some 
idempotent depending only on $\la$, see \cite[(7.7)]{J}. It follows that 
$$\SF(\si,\la) \otimes L(\tau,(n)) \simeq \SF(\si\tau,\la).$$ Each $\SF(\si,\la)$ has a unique maximal 
submodule, the quotient by which is exactly $L(\si,\la)$. Hence  
$L(\si,\la) \otimes L(\tau,(n)) \simeq L(\si\tau,\la)$. 
\end{proof}


\subsection{Lower bound for $\mathbf{\kappa_{\ell}}$}\label{SLB} Throughout $\S$\ref{SLB},  $\si$ denotes an element of $\bar\FF_q^\times$ of degree $d|n$ and $k=n/d$. 


\begin{Proposition}\label{red1}
 Let $m\in\ZZ_{\geq 1}$, $\si \in \overline{\FF}_{q}^{\times}$, and
$SL_{n} \leq R \leq GL_n$. If $\kappa_R(\SC(\si,\la))\geq m$ for all 
$\la \vdash k$, then $\kappa_R(L(\si,\la))\geq m$ for all $\la \vdash k$.
\end{Proposition}

\begin{proof}
By \cite[8.2]{J}, composition factors of $\SF(\si,\la)$ are of the form $L(\si,\mu)$ 
for $\mu \tri \la$, and exactly one of them is  
$L(\si,\la)$. We apply induction on the dominance order $\tri$ on partitions 
of $k$. If $\la=(k)$, then $L(\si,\la) = \SF(\si,\la)$, and 
the result follows, as $\SF(\si,\la)$ is a reduction modulo $\ell$ of $\SC(\si,\la)$. 

For the induction step, we assume that for each $\mu \tri \la$ and $\mu \neq \la$, we have 
$L(\si,\mu)\dar_{R} =  \oplus^{m(\mu)}_{i=1}L(\si,\mu,i)$ for some $m(\mu) \geq m$ and 
$L(\si,\mu,i) \in \IBr(R)$. Since $GL_n/R$ is cyclic, the inertia group of $L(\si,\mu,i)$ has index $m(\mu)$ in $GL_n$. By our assumption,
$\SC(\si,\la)\dar_{R} = \oplus^{t}_{j=1}\VC(j)$ for some 
integer $t \geq m$ and some irreducible $\CC R$-modules $\VC(j)$, $1 \leq j \leq t$, which are 
$GL_n$-conjugate. Let $V(j)$ denote an $\ell$-modular reduction of $\VC(j)$. Now in the Grothendieck 
group of $\FF GL_n$-modules we have
$$\SF(\si,\la) = L(\si,\la) + \sum_{\mu \tri \la,~ \mu \neq \la}d_{\la\mu}L(\si,\mu)$$
for some integers $d_{\la\mu} \geq 0$. Thus over $R$ we have
\begin{equation}\label{toR}
  \sum^{t}_{j=1}V(j) = L(\si,\la)\dar_{R} + \sum_{\mu \tri \la,~ \mu \neq \la}d_{\la\mu}
    \sum^{m(\mu)}_{i=1}L(\si,\mu,i).
\end{equation}

Assume for a contradiction that 
$L(\si,\la)\dar_{R} = L_{1} \oplus \ldots \oplus L_{s}$, where 
$L_{1}, \ldots ,L_{s}$ are distinct $GL_n$-conjugates of $L_{1} \in \IBr(R)$, and $s<m$. 
We may assume that $L_{1}$ is a composition factor of $V(1)$. 
For any $j > 1$, $\VC(j) = ~ ^{g}\VC(1)$ for some $g \in GL_n$,
and so $^{g}L_{1}$ is a subquotient of $V(j)$. Thus the left-hand side of (\ref{toR}) contains 
at least $t$ conjugates of $L_{1}$. Since $t \geq m > s$, (\ref{toR}) now 
implies that some $GL_n$-conjugate $^{g}L_1$ of $L_{1}$ must be isomorphic to some $L(\si,\mu,i)$. 
So the inertia group of $^{g}L_{1}$ has index $m(\mu) \geq m > s$ in $GL_n$,
a contradiction. 
\end{proof}

\begin{Corollary}\label{red2}
Let $\la\vdash k$, 
$1 \neq \al \in \FQ^{\times}$, and assume that $\si\al$ is conjugate to~$\si$. Then 
$\kappa_R(L(\si,\la))\geq |\al|$ for $R := \Ker(\SC(\al,(n)))$.
\end{Corollary}
\begin{proof}
By Lemma \ref{tensor}, we have 
$\SC(\si,\la) \otimes \SC(\al,(n)) \simeq \SC(\si\al,\la) \simeq \SC(\si,\la).$ 
Note also that $|\al| = (GL_n:R)$. By Lemma \ref{linear}(ii), 
$\kappa_{R}(\SC(\si,\la))\geq  |\al|$, and this is true for all 
$\la \vdash k$. Now the statement follows from Proposition 
\ref{red1}.    
\end{proof}

\begin{Lemma}\label{LNT}
Let $r,c,d$ be positive integers and  $\ell^c$ be the $\ell$-part of $r-1$. Then the 
$\ell$-part of $\frac{r^{\ell^d}-1}{r-1}$ is $\ell^d$, unless $c=1$, $\ell=2$ and 
$r\equiv 3\pmod{4}$.
\end{Lemma}
\begin{proof}
We only explain the case $c\geq 2$. As $\frac{r^{\ell^d}-1}{r-1}=\prod_{i=1}^{d}\frac{r^{\ell^i}-1}{r^{\ell^{i-1}}-1}$, it suffices to see that $\left(\frac{r^{\ell^i}-1}{r^{\ell^{i-1}}-1}\right)_\ell=\ell$ for $1\leq i\leq d$. For such $i$ we have $r^{\ell^{i-1}}=1+A\ell^b$ for some integer $b\geq 2$ and $(A,\ell)=1$. So $r^{\ell^{i}}=(1+A\ell^b)^\ell=1+A\ell^{b+1}+B\ell^{2b}$ for some $B\in\ZZ$, whence the $\ell$-part of $r^{\ell^i}-1$ is $\ell^{b+1}$ (using $2b>b+1$). It follows that  $\left(\frac{r^{\ell^i}-1}{r^{\ell^{i-1}}-1}\right)_\ell=\frac{\ell^{b+1}}{\ell^b}=\ell$, as required.  
\end{proof}

\begin{Lemma}\label{red3}
 Assume that $\si$ is an $\ell'$-element and  $\ell^{c}|\gcd(n,q-1)$ for some $c \geq 1$. Then for any $\la \vdash k$ with 
$\ell^{c}|\la'$, we have $\kappa_{R_n}(L(\si,\la))\geq\ell^{c}$.  
\end{Lemma}

\begin{proof}
As $\ell^{c}|\la'$, we have $\ell^c|k$, and so $\ell^{c} d|n$, where 
$d := \deg(\si)$. Let $P := \OLB(\FF_{q^{d\ell^{c}}}^{\times})$ and $Q := \OLB(\FF_{q^{d}}^{\times})$. Note that $Q<P$ via $\FF_{q^{d}}\subset \FF_{q^{d\ell^{c}}}$, and  $\ell^{c}$ divides $|P/Q|$. 

First, assume that $c = 1$. Then there is $\tau \in P \setminus Q$ such that 
$Q := \langle \tau^{\ell} \rangle$. Note that $\tau$ has degree $\ell$ over $\FF_{q^{d}}$, 
whence $\FF_{q^{d}}(\tau) = \FF_{q^{\ell d}}$. The $\ell'$-part (resp. $\ell$-part) of $\sigma\tau$ is $\si$ (resp. $\tau$). In particular, $\si$ and $\tau$ are powers of $\si\tau$. Therefore 
$$\FQ(\si\tau) = \FQ(\si,\tau) = \FQ(\si)(\tau) = \FF_{q^{d}}(\tau) = \FF_{q^{\ell d}},$$ 
i.e. $\deg(\si\tau) = \ell d$.
Set $\al := (\si\tau)^{q^{d}-1}$. We claim that $\al \in \FQ^{\times}$ and $|\al| = \ell$. Indeed, since 
$\si^{q^{d}-1} = 1$, we have $\al = \tau^{q^{d}-1}$. Now 
$\al^{\ell} = (\tau^{\ell})^{q^{d}-1} = 1$ since $\tau^{\ell} \in Q$;
in particular, $\al \in \FQ^{\times}$. Also, 
$\al \neq 1$, as otherwise $\tau \in \FF_{q^{d}}^{\times}$, a contradiction.

We have shown that $\si\tau$ has degree $\ell d$ and $(\si\tau)^{q^{d}} = \al(\si\tau)$ with 
$1 \neq \al \in \FQ^{\times}$. By Corollary \ref{red2}, $\kappa_K(L(\si\tau,\nu))\geq \ell$ for $K := \Ker(\SC(\al,(n))) \geq R_n$ and  
all $\nu \vdash n/\ell d$. In particular, $\kappa_K(L(\si\tau,\mu))\geq\ell$, where $\mu \vdash n/\ell d$ is chosen so that $\mu' = (1/\ell)\la'$. Since
$L(\si\tau,\mu) \simeq L(\si,\la)$ by \cite{DJ1}, the claim follows.     

Now let $c \geq 2$. Observe that $\ell^{c} = |P/Q|$ or  $\left(\frac{q^{d\ell^c}-1}{q^d-1}\right)_\ell=\ell^c$ by Lemma~\ref{LNT}. 
Now if 
$P=\langle\tau\rangle$, then $Q = \langle \tau^{\ell^{c}} \rangle$ and $\tau^{\ell^{c-1}} \notin Q$. Note $\tau$ has 
degree $\ell^{c}$ over $\FF_{q^{d}}$, 
so $\FF_{q^{d}}(\tau) = \FF_{q^{d\ell^c}}$. As the $\ell'$-part and the $\ell$-part of $\si\tau$ are $\si$ and $\tau$, respectively, $\si$ and $\tau$ are powers of $\si\tau$. So   
$$\FQ(\si\tau) = \FQ(\si,\tau) = \FQ(\si)(\tau) = \FF_{q^{d}}(\tau) = \FF_{q^{d\ell^{c}}},$$ i.e. $\deg(\si\tau) = d\ell^{c}$.
Set $\al := (\si\tau)^{q^{d}-1}=\tau^{q^{d}-1}$. We claim $\al \in \FQ^{\times}$ has order $\ell^{c}$. Indeed, 
$\al^{\ell^{c}} = (\tau^{\ell^{c}})^{q^{d}-1} = 1$ since 
$\tau^{\ell^{c}} \in Q$, and $\al \in \FQ^{\times}$ as $\ell^{c}|(q-1)$. On the other hand, 
$\al^{\ell^{c-1}} \neq 1$, as otherwise $\tau^{\ell^{c-1}} \in \FF_{q^{d}}^{\times}$ and so 
$\tau^{\ell^{c-1}} \in Q$, a contradiction.
Thus $\deg(\si\tau)=d\ell^{c}$ and $(\si\tau)^{q^{d}} = \al(\si\tau)$ for 
$\al \in \FQ^{\times}$ of order $\ell^{c}$. Now we finish as in the case $c=1$.     
\end{proof}

\begin{Theorem}\label{sl1}
Let $\ell^c |\gcd(n,q-1)$ and 
$V = L(\si_{1},\la^{(1)}) \circ \ldots L(\si_{a},\la^{(a)})$ be an irreducible $\FF GL_{n}$-module, where 
$\si_{1}, \ldots ,\si_{a}$ are $\ell'$-elements and $\ell^c|(\la^{(i)})'$ for all $i = 1,2, \dots ,a$. 
Then $\kappa_{R_n}(V)
\geq \ell^c$.
\end{Theorem}

\begin{proof}
We 
apply induction on $a$. The case $a = 1$ is Lemma~\ref{red3}. For the induction step, assume $a \geq 2$, and let $\la_{i} \vdash k_{i}, \deg(\si_{i})=d_{i}$. Set $r := k_{1}d_{1}$, $s := n-r,$ 
$A := GL_{r}$, $B := GL_{s}$, $A_{1} := R_{r}$, $B_{1} := R_{s}$, and
$$U := L(\si_{1},\la^{(1)}) \in \IBr(A),\quad W := (L(\si_{2},\la^{(2)}) \circ \ldots \circ L(\si_{a},\la^{(a)})) \in \IBr(B).$$ 
As $\ell^c|(\la^{(i)})'$ for all $i$, we have $\ell^c|\gcd(r,q-1)$ and $\ell^c|\gcd(s,q-1)$. 
By the inductive assumption, 
$$u := \kappa_{A_{1} \times B}(U \otimes W) = \ell^{\al} \geq \ell^c,$$ 
and 
$$v := \kappa_{A \times B_{1}}(U \otimes W) = \ell^{\beta} \geq \ell^c.$$ 
Choose $x \in GL_{r}$, 
$y \in GL_{s}$ such that $\langle \det(x) \rangle = \OLB(\FQ^{\times})$ and 
$\det(y) = \det(x)^{-1}$. Then $\bar{x} := xA_{1}$ generates $A/A_{1}$,
$\bar{y} := yB_{1}$ generates $B/B_{1}$, and $A/A_{1} \simeq B/B_{1}$. 

Next, consider the standard parabolic subgroup $P = QL<GL_n$ with Levi subgroup 
$L:= GL_{r} \times GL_{s}$. Note that $K := L \cap R_n=\langle A_{1},B_{1},xy \rangle$, i.e. $K/(A_{1} \times B_{1}) = \langle \bar{x}\bar{y} \rangle$.
By Lemma \ref{cross}, 
$$\kappa_{K}(U \otimes W)=\gcd(u,v) = \gcd(\ell^{\al},\ell^{\beta}) 
\geq \ell^c.
$$ 
Now $V = U \circ W$ and $(R_n:QK) = (GL_n:P)$, so 
$$V\dar_{R_n}\simeq (\infl_K^{QK}(U \otimes W))\uar^{R_n}.$$ Hence $\kappa_{R_n}(V)\geq \kappa_K(U \otimes W)\geq\ell^c$.   
\end{proof}

\begin{Theorem}\label{TRn}
Let $V = L(\sigma_{1},\la^{(1)}) \circ \ldots \circ L(\sigma_{a},\la^{(a)})$ be 
an irreducible $\FF GL_n$-representation, where the $\sigma_{i}$ are $\ell'$-elements. Set $\underline{\la}=(\la^{(1)}, \dots,\la^{(a)})$. Then $\kappa_{R_n}(V)=\gcd\big\{(GL_n:R_n),\Delta(\underline{\la}')\big\}$.
 \end{Theorem}
\begin{proof}
Note that $\gcd\big\{(GL_n:R_n),\Delta(\underline{\la}')\big\}$ is some power $\ell^{c(V)}$, which divides $\gcd(n,q-1)$, and all $(\la^{(i)})'$. By Theorem~\ref{sl1}, we have $\kappa_{R_n}(V)\geq \ell^{c(V)}.$

To see that $\kappa_{R_n}(V)= \ell^{c(V)}$ for each $V$, we apply a counting argument.  First apply Lemma~\ref{LGain} with $S=B=R_n$ and $G=A=GL_n$ to conclude that the set of all irreducible summands of the restriction $W\dar_{R_n}$, as $W$ runs over the  isomorphism classes $[W]$ of irreducible $\FF GL_n$-modules, is exactly the set of representatives of the isomorphism classes of the irreducible $\FF R_n$-modules. Hence it suffices to show that the number of $\ell$-regular conjugacy classes in $R_n$ equals $\sum_{[W]}\ell^{c(W)}$. 

In order to compute the number of $\ell$-regular conjugacy classes in $R_n$, we recall first that such classes in $GL_n$ are parametrized by the pairs $(s,\underline{\mu})$ where $s\in\mathcal{S}_{\ell'}$ and $\underline{\mu}\vdash\underline{k}(s)$. Moreover, 
by Proposition~\ref{PIndex}, such class splits into exactly $\gcd\big((GL_n:R_n),\Delta(\underline{\mu})\big)$ $R_n$-conjugacy classes, which is precisely $c(W)$ for the irreducible $\FF GL_n$-module $W=L(s,\underline{\mu}')$. We conclude that the total number of $\ell$-regular conjugacy classes in $R_n$ equals $\sum_{[W]}\ell^{c(W)}$.
\end{proof}

\subsection{Proof of Theorem~\ref{TMain}}
We first apply Lemma~\ref{rami} to the case $G=GL_n$ and $S=SL_n$. In particular $B=R_n$. We conclude that $\kappa_S(V)=\kappa_A(V)\cdot\kappa_B(V)$. But by Lemmas~\ref{tensor} and \ref{linear}(i), $\kappa_A(V)=\kappa_{\ell'}(\sfr)$. On the other hand, 
by Theorem~\ref{TRn}, $\kappa_B(V)=\kappa_{\ell}(\sfr)$. The proof is complete.

\section{An application: irreducible reductions modulo $\ell$}
Here we classify 
complex representations of $SL_n(q)$ 
irreducible modulo $\ell$. 

\subsection{A canonical composition factor} 
We will show that 
reduction modulo $\ell$ of 
$\SC(\sfr)$ has 
a canonical composition factor 
$L(\sfrs)$. 
The proof of the following lemma is an easy exercise.

\begin{Lemma}\label{obv}
Let $\al^{(i)},\beta^{(i)} \vdash n_{i}$ with $\alpha^{(i)} \tri \beta^{(i)}$ for $i = 1, \ldots ,m$. Then:
\begin{enumerate}
\item[{\rm (i)}] $\al^{(1)} + \ldots + \al^{(m)} \tri \beta^{(1)} + \ldots + \beta^{(m)}$;
\item[{\rm (ii)}] if $\al^{(1)} + \ldots + \al^{(m)} = \beta^{(1)} + \ldots + \beta^{(m)}$ then 
$\al^{(i)} = \beta^{(i)}$ for all $i$. 
\hfill $\Box$
\end{enumerate} 
\end{Lemma}

\begin{Lemma}\label{cf1}
Any 
composition factor of $\VC := \SC(\si,\al^{(1)}) \circ \ldots \circ \SC(\si,\al^{(m)})$ is of the form $\SC(\si,\la)$ with
$\la'\triop \sum^{m}_{i=1}(\al^{(i)})'$. Moreover, $\SC(\si,\la)$ is a multiplicity one 
composition factor of $\VC$ if\, $\sum^{m}_{i=1}(\al^{(i)})' = \la'$.
\end{Lemma}

\begin{proof}
By \cite{DJ2}, \cite{J2}, or \cite{BDK}, we have 
$$\SC(\si,\al) \circ \SC(\si,\beta) = \sum_{\gamma \vdash |\al|+|\beta|}a_{\al\beta\gamma}\SC(\si,\gamma),$$
where $a_{\al\beta\gamma}$ are Littlewood-Richardson coefficients. Now the result follows by induction on $m$ using associativity of Harish-Chandra induction and well-known properties of Littlewood-Richardson coefficients. 
\end{proof}

\begin{Lemma}\label{cf2}
If $\si$ be an $\ell'$-element, then every composition factor of  
$V: = L(\si,\la^{(1)}) \circ \ldots \circ L(\si,\la^{(m)})$ 
is $L(\si,\ga)$, with multiplicity one, or $L(\si,\eta)$ with $\eta \tri \ga$ and 
$\eta \neq \ga$, where  
$\sum^{m}_{i=1}(\la^{(i)})' = \ga'$. 
\end{Lemma}

\begin{proof}
Let $k=\deg(\si)$. The unitriangularity of the submatrix $\Delta(\si,k)$ (see \cite{DJ2}) 
of the decomposition matrix of $GL_{kd}$ implies that in the Grothendieck
group of $\FF GL_{kd}$-modules we have 
$$L(\si,\la^{(i)}) = \overline{\SC(\si,\la^{(i)})} + \sum_{\al \tri \la^{(i)},~ \al \neq \la^{(i)}} 
  x_{i\al}\overline{\SC(\si,\al)}\qquad (x_{i\al}
  \in \ZZ).$$
Hence
\begin{equation}\label{forV}
  V = \overline{V}_{\CC} + \sum
    y_{\al^{(1)},...,\al^{(m)}}\overline{\SC(\si,\al^{(1)}) \circ \ldots
      \circ \SC(\si,\al^{(m)})},
\end{equation}
where the sum is over all $\al^{(i)} \tri \la^{(i)}$ such that $\al^{(j)} \neq \la^{(j)}$ for at least one $j$,  
$\VC = \SC(\si,\la^{(1)}) \circ \ldots \circ \SC(\si,\la^{(m)})$ and 
$y_{\al^{(1)},...,\al^{(m)}} \in \ZZ$.

Assume that $L(\si,\eta)$ is a composition factor of some summand
$$\overline{\SC(\si,\al^{(1)}) \circ \ldots \circ \SC(\si,\al^{(m)})}$$
in the summation in (\ref{forV}). 
Then by Lemma \ref{cf1},  
$L(\si,\eta)$ is a composition factor 
of some $\overline{\SC(\si,\beta)}$, where 
$\sum^{m}_{i=1}(\al^{(i)})' \tri \beta'$. Then $\eta \tri \beta$, whence 
$\beta' \tri \eta'$. But $(\la^{(i)})' \tri (\al^{(i)})'$ for each 
$i$ (since $\al^{(i)} \tri \la^{(i)}$). So by Lemma \ref{obv}(i),
$$\sum^{m}_{i=1}(\la^{(i)})' ~\tri~ \sum^{m}_{i=1}(\al^{(i)})' ~\tri~ \beta' 
  ~\tri~ \eta'.$$ 
Since $\sum^{m}_{i=1}(\la^{(i)})' = \ga'$ by the choice of $\ga$,
we get $\eta \tri \ga$. In particular, if $\eta = \ga$, then  
 Lemma \ref{obv}(ii)
implies that $(\la^{(i)})' = (\al^{(i)})'$, and so $\la^{(i)} = \al^{(i)}$ for 
all $i$, a contradiction. 

On the other hand, by Lemma \ref{cf1}, $\SC(\si,\ga)$ is a multiplicity one 
composition factor of $\VC$, with other composition factors of the form 
$\SC(\si,\mu)$ for $\mu\tri\ga$, 
$\mu\neq\ga$. Finally, $\overline{\SC(\si,\mu)}$ contains $L(\si,\mu)$ as a 
multiplicity one composition factor, with other factors of the form $L(\si,\eta)$ 
for $\eta\tri\mu$, $\eta\neq\mu$. 
So the right hand of (\ref{forV}) contains $L(\si,\ga)$ with multiplicity one,
and all other composition factors are of
the form $L(\si,\eta)$ with $\eta \tri \ga$ and $\eta \neq \ga$,
as stated.
\end{proof}

\begin{Theorem}\label{cf-main}
Let $V=L_\CC(\sfr)$ be an irreducible $\CC GL_{n}$-module.
Then $V(\mod\ell)$ has $L(\sfrs)$ as a composition factor with multiplicity one.
\end{Theorem}

\begin{proof}
In the notation of \S1.2, assume that the $\ell'$-parts $s_{i}$ of $\si_{i}$ are conjugate to
$s_{1}$ precisely when $i\in \{ j_{1} = 1, j_{2}, \ldots ,j_{b}\}$. 
By \cite{DJ2}, the Brauer character of 
$$\SC(\si_{j_{1}},\la^{(j_{1})}) \circ \ldots \circ \SC(\si_{j_{b}},\la^{(j_{b})})(\mod\ell)$$  
is a {\em sum} of the Brauer character of  
$W_{1} := L(s_{1},\mu^{(j_{1})}) \circ \ldots \circ L(s_{1},\mu^{(j_{b})})$ 
with some other Brauer character $\psi$, where 
$(\mu^{(i)})' = k_{i}(\la^{(i)})'$. By Lemma \ref{cf2}, $W_{1}$ has $L(s_{1},\delta^{(1)})$ as a multiplicity one composition factor (see 
\S1.2). On the other hand, using~\cite{DJ1},  
\begin{eqnarray*}
\overline{\SC(\si_{j_{i}},\la^{(j_{i})})}&=& L(\si_{j_{i}},\la^{(j_{i})}) 
  + \sum_{\beta^{(j_{i})} \tri \la^{(j_{i})},~\beta^{(j_{i})} \neq \la^{(j_{i})}}
    x^{i}_{\beta^{(j_{i})}}L(\si_{j_{i}},\beta^{(j_{i})}) \\&=&
L(s_{1},\mu^{(j_{i})}) 
  + \sum_{\beta^{(j_{i})} \tri \la^{(j_{i})},~\beta^{(j_{i})} \neq \la^{(j_{i})}}
    x^{i}_{\beta^{(j_{i})}}L(s_{1},\nu^{(j_{i})}),
\end{eqnarray*}
where $x^{i}_{\beta^{(j_{i})}} \in \ZZ$ and 
$(\nu^{(j_{i})})' = k_{j_{i}}(\beta^{(j_{i})})' \triop (\mu^{(j_{i})})'$.
Hence, by Lemma \ref{cf2}, 
any composition factor of $\psi$ is of the form 
$L(s_{1},\gamma)$, where 
$\gamma \tri \delta^{(1)}$ but $\gamma \neq \delta^{(1)}$.  
Now the result comes from 
repeating this argument for all other conjugate classes of $s_{i}$.
\end{proof}

\subsection{The James-Mathas Theorem}
$\CC GL_{n}$-modules irreducible modulo $\ell$ have been classified
by James and Mathas. 
For reader's convenience, we reformulate the 
result here adapting to the notation of \S1.2.   

\begin{Theorem} {\rm \cite[Theorem (4.33)]{JM}}\label{jm}
The reduction $\SC(\sfr)(\mod \ell)$ is irreducible
if and only the following two conditions hold:

{\rm (JM1)} If $i \neq j$ and $\deg(\si_{i}) = \deg(\si_{j})$, then 
$[s_{i}] \neq [s_{j}]$;

{\rm (JM2)} For all $i$ and all nodes $(a,c)$ and $(b,c)$ in the Young diagram of
$\la^{(i)}$, the corresponding hook lengths $h_{ac}$ and $h_{bc}$ satisfy the condition
$$\left(\frac{q^{h_{ac}f}-1}{q^{f}-1}\right)_{\ell} = 
  \left(\frac{q^{h_{bc}f}-1}{q^{f}-1}\right)_{\ell}$$
where $f := \deg(\si_{i})$.  
\end{Theorem} 

A partition $\la$ will be called a {\it $JM$-partition}
if it satisfies the condition (JM2) for some power $q^{f} > 1$. (Empty partition is a JM-partition).

\begin{Lemma}\label{jm-p}  
Let $\ell|(q-1)$ and $\la$, $\la^{(1)},\dots,\la^{(m)}$, $\mu$, $\mu^{(1)},\dots,\mu^{(n)}$ be $JM$-partitions, with all the $\la^{(i)},\mu^{(j)}$ non-empty. Then:
\begin{enumerate}
\item[{\rm (i)}] if $\la' \equiv \mu' (\mod \ell)$, then $\la = \mu$;
\item[{\rm (ii)}] if $\ell {{|}} \la'$, then $\la$ is empty;
\item[{\rm (iii)}] if $0 \leq a_{1} < a_{2} < \ldots < a_{m}$,  
$0 \leq b_{1} < b_{2} < \ldots < b_{n}$ are integers and 
$\sum^{m}_{i=1}\ell^{a_{i}}(\la^{(i)})' = 
  \sum^{n}_{j=1}\ell^{b_{j}}(\mu^{(j)})'$, 
then $m = n$ and $a_{i} = b_{i}$, $\la^{(i)} = \mu^{(i)}$ for all~$i$.  
\end{enumerate}
\end{Lemma}

\begin{proof}
(i) Adding 
empty columns to 
$\la$ and $\mu$, we
may assume that both $\la$ and $\mu$ have 
$c$ columns and the $c$th column of both is empty. We prove by induction
on $i=0,1,\dots, c$, that $\la'_{c-i} = \mu'_{c-i}$. The induction base is clear. 
So suppose 
$\la'_{j} = \mu'_{j}$ for $c-k+1 \leq j \leq c$, 
and prove 
that $\la'_{c-k} = \mu'_{c-k}$. By assumption $\la'_{c-k} \equiv \mu'_{c-k} (\mod \ell)$, and 
we may assume $\la'_{c-k} \geq \mu'_{c-k}$. 
So $\la'_{c-k} \geq \mu'_{c-k} + \ell$. Then we have
$$\la'_{c-k+1} = \mu'_{c-k+1} \leq \mu'_{c-k} \leq \la'_{c-k}-\ell.$$ Hence, setting 
$b := \la'_{c-k}$ and $a := b-\ell+1$, we have $h_{b,c-k} = 1$ and $h_{a,c-k+1} = \ell$.
Now, $\la$ satisfies the condition 
(JM2) for some $Q := q^{f}>1$. As  $\ell|(q-1)$, we have $\ell|(Q-1)$, so $(Q^{h_{b,c-k}}-1)/(Q-1) = 1$, but $(Q^{h_{a,c-k}}-1 )/(Q-1) = (Q^{\ell}-1)/(Q-1)$ is divisible by $\ell$, contradicting (JM2). 
   
(ii) follows from (i) by taking $\mu = \emptyset$.

(iii) Adding 
empty partitions, we may assume that $m = n$ and $a_{i} = b_{i} = i-1$.
The 
equality 
$\sum^{n}_{i=1}\ell^{i-1}(\la^{(i)})' = \sum^{n}_{i=1}\ell^{i-1}(\mu^{(i)})'$ now implies 
$(\la^{(1)})' \equiv (\mu^{(1)})' (\mod \ell)$. By (i), we get $\la^{(1)} = \mu^{(1)}$. 
Now we have 
$\sum^{n-1}_{i=1}\ell^{i-1}(\la^{(i+1)})' = \sum^{n-1}_{i=1}\ell^{i-1}(\mu^{(i+1)})'$,
and the claim follows by induction on $n$.  
\end{proof}

\subsection{Relations between branching numbers}
\begin{Lemma}\label{parts}
Let $\si, \si' \in \overline{\FF}_{q}^{\times}$, $\si_\ell=u$,
$\si_{\ell'}=s$, $\si'_{\ell}=u'$, $\si'_{\ell'}=s'$, and 
$[\si \tau] = [\si']$ for some $\tau \in \FQ^{\times}$. Then $\deg(\si) = \deg(\si')$. Moreover:
\begin{enumerate}
\item[{\rm (i)}] if $\tau$ is an $\ell$-element then $[s] = [s']$;
\item[{\rm (ii)}] if $\tau$ is an $\ell'$-element then $[s\tau] = [s']$ and 
$\deg(s) = \deg(s')$.
\end{enumerate}
\end{Lemma}

\begin{proof}
The first statement is clear. By assumption, $\si \tau = (\si')^{q^{i}}$ for some $i \geq 0$. Then in (i) we have that 
$\tau u(u')^{-q^{i}} = (s')^{q^{i}}s^{-1}$ is both 
an $\ell$-element and an $\ell'$-element and so it must be $1$. Hence  
$s = (s')^{q^{i}}$, i.e. $[s] = [s']$. 
In (ii) we have that 
$u(u')^{-q^{i}} = (s')^{q^{i}}s^{-1}\tau^{-1}$ is both 
an $\ell$-element and an $\ell'$-element and so it must be $1$. It follows that 
$s\tau = (s')^{q^{i}}$, i.e. $[s\tau] = [s']$. 
\end{proof}

\begin{Lemma}\label{k2}
If $\si \in \overline{\FF}_{q}^{\times}$, $s=\si_{\ell'}$, $\deg(\si)=kd$, and $\deg(s)=d$, then the order of the group 
$I = \{\tau \in \OLB(\FQT ) \mid 
  [\si]=[\si \tau]\}$ divides $\gcd(k,q-1)_{\ell}$. In fact, $|I| = \gcd(k,q-1)_{\ell}$, except
for the case where $\ell = 2$, $q^{d} \equiv 3 (\mod 4)$, $k = 2$, and the 
$\ell$-part $u$ of $\si$ has order $\geq 8$. In the exceptional case $|I| = 1$.
\end{Lemma}

\begin{proof}
We may assume that $\ell|(q-1)$. 
For $\tau \in I$ we have $su\tau = \si \tau = (\si)^{q^{j}} = s^{q^{j}}u^{q^{j}}$
for some $j \geq 0$. So 
$s^{1-q^{j}} = u^{q^{j}-1}\tau^{-1}=1$. As $\deg(s) = d$, we must have 
$d|j$. Hence, using the facts that $\tau^q=\tau$ and $\deg(\si) = kd$, we get
$$\tau^{k} = \tau^{1+q^{j} + \ldots + q^{(k-1)j}} = 
  \left(\si^{q^{j}-1}\right)^{1+q^{j} + \ldots + q^{(k-1)j}} = 
  \si^{q^{kj}-1} = 1.$$
Now $|I|$ divides $N := \gcd(k,q-1)_{\ell}$, as $I$ is a cyclic $\ell$-subgroup of $\FQT $. 

The second statement is now obvious if $k = 1$, so we may assume by \cite{DJ1} that 
$k = \ell^{a}$ and $N = \ell^{b}$ for some integers $a \geq b \geq 1$.  
Suppose first that $\ell \neq 2$ if $q^{d} \equiv 3 (\mod 4)$. Then by Lemma~\ref{LNT}, for any integers $c,m \geq 1$, 
$(q^{md\ell^{c}}-1)/(q^{md}-1)$ has the $\ell$-part equal to $\ell^{c}$. 
Set $$\tau = \si^{q^{d\ell^{a-b}}-1} = u^{q^{d\ell^{a-b}}-1}.$$
In particular, $\tau$ is an $\ell$-element.
Since $\deg(\si) = d\ell^{a}$ and the $\ell$-part of 
$(q^{d\ell^{a}}-1)/(q^{d\ell^{a-b}}-1)$ is $\ell^{b}$, we see that 
$\tau^{\ell^{b}} = 1$. But $\ell^{b} | (q-1)$, hence $\tau^{q-1} = 1$ and so $\tau \in I$.
On the other hand, if $\tau^{\ell^{b-1}} = 1$, then, since the $\ell$-part of 
$(q^{d\ell^{a-1}}-1)/(q^{d\ell^{a-b}}-1)$ is $\ell^{b-1}$, we must have  
$$\si^{q^{d\ell^{a-1}}-1} = \tau^{(q^{d\ell^{a-1}}-1)/(q^{d\ell^{a-b}}-1)} =1,$$
contrary to the equality $\deg(\si) = \ell^{a}d$. Thus $|\tau| = \ell^{b} = N$ and
so $|I| = N$.   

Now let $\ell = 2$ and $q^{d} \equiv 3 (\mod 4)$; in particular,
$N = 2$. First we consider the case where either $a \geq 2$, or $a = 1$ but 
$|u| < 8$ (hence $|u| = 4$). Set 
$$\tau = \si^{q^{2^{a-1}d}-1} = u^{q^{2^{a-1}d}-1}.$$
Then $\tau \neq 1$ since $\deg(\si) = 2^{a}d$. On the other hand, if $a \geq 2$, 
then the $2$-part of $(q^{2^{a}d}-1)/(q^{2^{a-1}d}-1)$ is $2$ and so 
$\tau^{2} = 1$, whence $\tau \in I$ and $|I| = N$. If $a = 1$ and $|u| = 4$, then
$\tau^{2} = u^{2(q^{d}-1)} = 1$, and again $|I| = N$. 

Finally, if  $a = 1$ and 
$I \ni \tau \neq 1$, then $\tau = -1$ and $-\si = (\si)^{q^{j}}$ for some
$j$, $1 \leq j < 2d$. As above, $d|j$, so $j = d$. Hence
$-1 = (\si)^{q^{d}-1} = u^{q^{d}-1}$. Since $q^{d} \equiv 3 (\mod 4)$ and 
$u$ is a $2$-element, we conclude that $|u| = 4$.       
\end{proof}

Let $V_\CC= \SC(\sfr)$ be an arbitrary irreducible $\CC GL_n$-module. By Theorem~\ref{TMain}, $\kappa^{GL_n}_{SL_n}(V_\CC)= \{\tau \in \FQT  \mid 
\tau\cdot\sfr=\sfr\}$. Denote the $\ell$-part (resp. $\ell'$-part) of this number by  $\kappa_\ell(V_\CC)$ (resp. $\kappa_{\ell'}(V_\CC)$). Thus $\kappa^{GL_n}_{SL_n}(V_\CC)=\kappa_\ell(V_\CC)\kappa_{\ell'}(V_\CC)$. 

\begin{Proposition}\label{k2-main}
Let $\VC = \SC(\sfr)$. 
Then:
\begin{enumerate}
\item $\kappa_\ell(V_\CC)$ 
divides $\kappa_\ell(\sfrs)$.
\item Assume in addition that $\VC (\mod \ell)$ is irreducible. If $\sfr$ is not critical, then 
$\kappa_\ell(V_\CC) = \kappa_\ell(\sfrs)$. If $\sfr$ is critical, then 
$\kappa_\ell(V_\CC) = 1$ and $\kappa_{\ell}(\sfrs) = 2$.
\end{enumerate}
\end{Proposition}

\begin{proof}
(i) Let $\tau$ be a generator of the group 
$I := \{\si \in \OLB(\FQT ) \mid 
\si\cdot\sfr=\sfr\}$. Then $\kappa_{\ell}(V_\CC)=|\tau| = \ell^{e}$ for some $e$. 
We adopt the notation $\si_{i} = s_{i}u_{i}$, $k_{i}$,
$d_{i}$, as in \S1.2. By Lemma \ref{parts}(i), 
if $[\si_{i}\tau] = [\si_{j}]$ 
then $[s_{i}] = [s_{j}]$. 
So the multiplication 
by $\tau$ preserves the set $\SCL_{1}$ of 
all $(\si_{i},\la^{(i)})$ with $[s_{i}] = [s_{1}]$. 
Consider one $\tau$-orbit $\OL\subset\SCL_{1}$, say $\tau$ permutes 
$$(s_{1}u_{1},\la^{(1)}),(s_{1}u_{2},\la^{(1)}),\ldots ,
(s_{1}u_{b},\la^{(1)})$$ 
cyclically (after a suitable relabeling of the 
$(\si_{i},\la^{(i)})$), with $b = \ell^{c}$ for some $0 \leq c \leq e$.
Then $\tau^{b}$ preserves $(s_{1}u_{1},\la^{(1)})$. By Lemma \ref{k2},
$|\tau^{b}| = \ell^{e-c}$ divides $\gcd(k_{1},q-1)_{\ell}$. Defining 
$\gamma_{\OL} := (bk_{1}(\la^{(1)})')'$, we see that $\ell^{e} = |\tau|$ divides 
$(\gamma_{\OL})'$. Notice that, in the definition of $\sfrs$, 
$(\delta^{(1)})'$ is just the sum of all the partitions $(\gamma_{\OL})'$, when 
$\OL$ runs over all orbits inside $\SCL_{1}$. Hence 
$\ell^{e}|(\delta^{(1)})'$. 
Similarly $\ell^{e}|(\delta^{(j)})'$ for all $j$,
and of course $\ell^{e}$ divides $n$ and $q-1$. So $\ell^{e}$ divides $\kappa_{\ell}(\sfrs)$.      

(ii) is clear if $\ell{\not{|}}(q-1)$, so let $\ell|(q-1)$. 
We apply~Theorem~\ref{jm}
to~$\VC$. Note that 
each orbit $\OL$ as above has size $1$, for, if $[\si_{i}\tau] = [\si_{j}]$, we saw that $\deg(\si_{i}) = \deg(\si_{j})$, and 
$[s_{i}] = [s_{j}]$, whence $i=j$ by 
(JM1). So the multiplication by $\tau$ 
preserves each $([\si_{i}],\la^{(i)})$. Moreover, if $[s_{i}] = [s_{j}]$ for some 
$i \neq j$, then (JM1) implies that $k_{i}d_{i} \neq k_{j}d_{j}$,
and so $k_{i} \neq k_{j}$ as $d_{i} = d_{j}$. So, if 
$$\SCL_{1} = \{(\si_{j_{1}},\la^{(j_{1})}), \ldots ,(\si_{j_{b}},\la^{(j_{b})})\},$$ 
then we may assume that 
$
k_{j_{1}} < k_{j_{2}} < \ldots < k_{j_{b}}$. On the other hand, 
by Lemma \ref{jm-p}(ii), 
we have $\ell{\not{|}}(\la^{(i)})'$ for all $i$. As the $k_{i}$ are $\ell$-powers, the exact $\ell$-power dividing    
$(\delta^{(1)})' = \sum^{b}_{i=1}k_{j_{i}}(\la^{(j_{i})})'$ is 
$\gcd(k_{j_{1}}, \ldots ,k_{j_{b}})$. A similar result for all $\de^{(j)}$ yields
\begin{equation}\label{fork2}
  \kappa_{\ell}(\sfrs) = \gcd(q-1,k_{1}, \ldots ,k_{a})_{\ell}.
\end{equation}

Let $J_{i}:=\{\si \in \OLB(\FQT )\mid [\si\si_i]=[\si_i]\}$. 
As noted above, $\tau\in J_i$ for each $i$, so $I = \cap^{a}_{i=1}J_{i}$.
Assume $\ell \neq 2$ if $q \equiv 3 (\mod 4)$. By Lemma \ref{k2}, 
$J_{i}$ 
is cyclic of order $\gcd(q-1,k_{i})_{\ell}$, so 
$|I| = \gcd(q-1,k_{1}, \ldots ,k_{a})_{\ell}$, as required.
Finally, let $\ell = 2$ and $q \equiv 3 (\mod 4)$. Then 
$|I| \leq \kappa_{\ell}(\sfrs) \leq 2$. So  
$|I| \neq \kappa_{\ell}(\sfrs)$ if and only if $|I| = 1$ and 
$\kappa_{\ell}(\sfrs) = 2$. By (\ref{fork2}), $\kappa_{\ell}(\sfrs) = 2$
if and only if all $k_{i} \geq 2$. On the other hand, since $I = \cap^{a}_{i=1}J_{i}$, and
each $J_{i}$ is a subgroup of the cyclic group $O_{2}(\FQ^{\times})$, we have 
$|I| = 1$ exactly when $|J_{j}| = 1$ for some $j$. As $2|k_{j}$, Lemma \ref{k2} implies that $|J_{j}| = 1$ precisely when $k_{j} = 2$, $d_{j}$ is odd,
and $|u_{j}| \geq 8$. We have shown that $|I| \neq \kappa_{\ell}(\sfrs)$ precisely when
$\sfr$ is critical.       
\end{proof}

\begin{Proposition}\label{k1-main}
Let $\VC = \SC(\sfr)$. Then $\kappa_{\ell'}(\VC)$ 
divides $\kappa_{\ell'}(\sfrs)$.
\end{Proposition}

\begin{proof}
Let 
$J := \{\nu \in \OLA(\FQT ) \mid 
  \nu\cdot\sfr=\sfr\}=\langle\rho\rangle$. 
 We adopt the notation 
 of \S1.2. By Lemma \ref{parts}(ii), 
if $[\si_{i}\rho] = [\si_{j}]$, then $[s_{i}\rho] = [s_{j}]$ and $k_{i} = k_{j}$.
So the multiplication 
by $\rho$ sends the set $\SCL_{1}=\{(\si_{i},\la^{(i)})\mid [s_{i}] = [s_{1}]\}$ to some $\SCL_{j}=\{(\si_{i},\la^{(i)})\mid [s_{i}] = [s_{j}]\}$, where $[s_j]= [s_{1}\rho]$. 
If $$\SCL_{1} = \{(\si_{i_{1}},\la^{(i_{1})}), \ldots ,(\si_{i_{b}},\la^{(i_{b})})\},$$ then 
$$\SCL_{j} =   \{(\rho\si_{i_{1}},\la^{(i_{1})}), \ldots ,(\rho\si_{i_{b}},\la^{(i_{b})})\},$$
and 
$$(\delta^{(1)})' = \sum^{b}_{v=1}k_{i_{v}}(\la^{(i_{v})})',~~~
   (\delta^{(j)})' = \sum^{b}_{v=1}k_{i_{v}}(\la^{(i_{v})})'.$$
So 
$\rho$ sends the component
$(s_{1},\delta^{(1)})$ of $\sfrs$ to the component $(s_{j},\delta^{(j)})$
of $\sfrs$. Thus $J$ stabilizes $\sfrs$, and so 
$\kappa_{\ell'}(\VC) = |J|$ divides $\kappa_{\ell'}(\sfrs)$.  
\end{proof}

\subsection{Proof of Theorem \ref{main2}}
Set $G = GL_{n}$, $S = SL_{n}$, and  
$\VC \dar_{S} = \sum^{t}_{i=1}W_\CC^{i}$, a 
sum of 
$t = \kappa_{\ell'}(\VC) \cdot \kappa_{\ell}(\VC)$ 
irreducibles. By 
Theorem~\ref{cf-main}, $V := L(\sfrs)$ is a composition factor of 
$\VC (\mod \ell)$. By Theorem \ref{TMain}, $V\dar_{S}$ is a sum of 
$t' := \kappa_{\ell'}(\sfrs) \cdot \kappa_{\ell}(\sfrs)$ irreducibles.
By Propositions \ref{k1-main}, 
\ref{k2-main}(i), $t|t'$. 
 
If $\WN(\mod \ell)$ is irreducible then $W_\CC^{i} (\mod \ell)$ 
are irreducible for all $i$. So $(\VC (\mod \ell))\dar_{S}$ has 
exactly $t$ composition factors, and $t \geq t'$. As $t|t'$, we have
$t = t'$ and $\VC (\mod \ell) = V$ is irreducible. Now $t = t'$ implies $\kappa_{\ell'}(V_\CC) = \kappa_{\ell'}(\sfrs)$ and 
$\kappa_{\ell}(V_\CC) = \kappa_{\ell}(\sfrs)$, 
whence $\sfr$ is not critical by Proposition
\ref{k2-main}(ii). 

Conversely, if $\VC (\mod \ell)$ is irreducible, 
$\kappa_{\ell'}(\sfr) = \kappa_{\ell'}(\sfrs)$, and $\sfr$ is not critical, then 
$\kappa_{\ell}(\sfr) = \kappa_{\ell}(\sfrs)$ by Proposition \ref{k2-main}(ii). 
So $\VC (\mod \ell) = V$, and $t = t'$. 
Hence 
$\WN (\mod \ell)$ is irreducible. 

\begin{Example}\label{E1}
{\rm (i) Let $\VC = \SC(\sfr)$,  
 $\WN$ be an irreducible constituent of $\VC\dar_{SL_{n}}$, and  
all $\si_{i}$ be $\ell'$-elements. Then 
{\it $\VC(\mod \ell)$ is irreducible over $GL_{n}$ if and only if 
$\WN (\mod \ell)$ is irreducible over $SL_{n}$}. Indeed, in this case 
$\sfr = \sfrs$ and that $\sfr$ is not critical.
So the claim follows fromTheorem~\ref{main2}.

\smallskip
(ii) The conditions (ii), (iii) of Theorem 
\ref{main2} cannot be relaxed. 
In the examples below $\WN = \VC \dar_{SL_{n}}$ is irreducible,
but $\WN(\mod\ell)$ is reducible.

(a) Let $q \equiv 3 (\mod 4)$, $\ell = 2$, $2|n$, and $d$ be an odd divisor of 
$n/2$. Pick $s \in \overline{\FF}_{q}^{\times}$ to be of order $(q^{d}-1)/2$. Then $\deg(s) = d$. Also, choose $u \in \FF_{q^{2d}}^{\times}$ 
of order $\geq 8$. Set $\sfr = [([su], (n/2d))]$. Then 
$\SC(\sfr)(\mod \ell)$ is irreducible, 
$\kappa_{\ell'}(\sfr) = \kappa_{\ell'}(\sfrs) = 1$ (as can be seen 
by direct computation), but $\sfr$ is critical.

(b) Let $r < \ell$ be primes with  $\ell r|(q-1)$. Pick an element 
$\eps \in \FQ^{\times}$ of order $r$, and $r$ distinct $\ell$-elements 
$1 = t_{1}, t_{2}, \ldots ,t_{r} \in \FQ^{\times}$. Define 
$$\sfr = [([\eps t_{1}],(n/r)),([\eps^{2} t_{2}],(n/r)), \ldots ,
  ([\eps^{r} t_{r}],(n/r))].$$
Then $\SC(\sfr)(\mod \ell)$ is irreducible, $\sfr$ is not critical,
but $\kappa_{\ell'}(\sfr) = 1$ (by Lemma \ref{parts}) and 
$\kappa_{\ell'}(\sfrs) = r$. 
}   
\end{Example}

\section{A unitriangular decomposition submatrix}\label{SDec}

We will need the following partial converse to Lemma \ref{k2}:
 
\begin{Lemma}\label{k2c}
Let $d \in \NN$ and $\ell^{a}|(q-1)$ for some $a\in \ZZ_{\geq 0}$. Then there 
exists an $\ell$-element $u = u(a,d) 
\in \overline{\FF}_{q}^{\times}$, depending only on $a,d$,  
such that, for any 
$\ell'$-element $s \in \overline{\FF}_{q}^{\times}$ of degree $d$, we have $\deg(su) = \ell^{a}d$ 
and the order of 
$I := \{ t \in \OLB(\FQT ) \mid [su] = [sut]\}$ is $\ell^{a}$.  
\end{Lemma}

\begin{proof}
Taking $u = 1$ in the case $a = 0$, 
we may assume that $a \geq 1$.
Let $P := \OLB(\FF_{q^{d\ell^{a}}}^{\times})$ and $Q := \OLB(\FF_{q^{d}}^{\times})$. 
Note that $Q<P$ via $\FF_{q^{d}}\subset \FF_{q^{d\ell^{a}}}$, $P$ is cyclic, and 
$\ell^{a}$ divides $|P/Q|$. We will distinguish two cases.

Case 1: Either $\ell \neq 2$, or 
$\ell = 2$ but $q^{d} \equiv 1 (\mod 4)$. Then $|P/Q| = \ell^{a}$ by Lemma \ref{LNT}, and we 
choose $u$ to be a generator of $P$. 

Case 2: $\ell = 2$ and $q^{d} \equiv 3 (\mod 4)$. 
Then $Q = \langle -1 \rangle \simeq C_{2}$, $\ell^{a} = 2$, and we choose $u \in P$ such that 
$u^{2} = -1$. 

In either case, we have  $u \in P$ with  
$u^{\ell^{a}} \in Q$ but $u^{\ell^{a-1}} \notin Q$. In particular, 
\begin{equation}\label{foru1} 
  u^{(q^{d}-1)\ell^{a}} = 1,\mbox{ but }u^{(q^{d}-1)\ell^{a-1}} \neq 1.
\end{equation}

We 
show that $\deg(su) = d\ell^{a}$ for any $s$ as in the assumption.
Note that $e:=\deg(su)$ is divisible by 
$d = \deg(s)$, as $s$ is the $\ell'$-part of $su$. Write $e = kd$ for some $k \in \NN$.
Then $k|\ell^{a}$ since $s,u\in \FF_{q^{d\ell^{a}}}$. If 
$k < \ell^{a}$ then $k|\ell^{a-1}$, so
\begin{equation}\label{foru2} 
  1 = (su)^{q^{d\ell^{a-1}}-1} = u^{q^{d\ell^{a-1}}-1}.
\end{equation}
In Case 1, the $\ell$-part of $\frac{q^{d\ell^{a-1}}-1}{q^{d}-1}$ is $\ell^{a-1}$ by Lemma 
\ref{LNT}. As $u$ is an $\ell$-element, (\ref{foru2}) implies 
$u^{(q^{d}-1)\ell^{a-1}} = 1$, contrary to (\ref{foru1}). In Case 2, 
we have $a = 1$, so 
(\ref{foru2}) yields $1 = u^{q^{d}-1}$, 
again contradicting
(\ref{foru1}). Thus $k = \ell^{a}$ as needed.
 
Set $t := (su)^{q^{d}-1} = u^{q^{d}-1}$. Then  $sut = (su)^{q^{d}}$
and so $[sut] = [su]$.  By (\ref{foru1}), $t^{\ell^{a}} = u^{(q^{d}-1)\ell^{a}} = 1$. 
But $\ell^{a}|(q-1)$, whence $t \in \OLB(\FQT )$ and so $t \in I$. Moreover,
$t^{\ell^{a-1}} = u^{(q^{d}-1)\ell^{a-1}} \neq 1$ by (\ref{foru1}). Thus  
$|t| = \ell^{a}$. On the other hand, $|I|$ divides $\gcd(\ell^{a},q-1)_{\ell} = \ell^{a}$ by
Lemma \ref{k2}. So $|I| = \ell^{a}$.  
\end{proof}

Let ${\mathfrak L}/O_{\ell'}(\FF_q)$ be the set of $\OLA(\FQT)$-orbits on $\LF$ 
or a set of the orbit representatives depending on the context. We define a relation 
$\tri$ on ${\mathfrak L}/O_{\ell'}(\FF_q)$. First, for 
$\sfr = [(\si_{1},\la^{(1)}), \ldots ,(\si_{a},\la^{(a)})]$ and 
$\tfr = [(\tau_{1},\mu^{(1)}), \ldots ,(\tau_{b},\mu^{(b)})]$ in $\LF$, 
we set $\sfr \tri \tfr$ if $a = b$, and, after a suitable relabeling of the 
$(\si_{i},\la^{(i)})$, there is $\nu \in \OLA(\FQT)$ such that 
$[\si_{i}\nu] = [\tau_{i}]$ and $\la^{(i)} \tri \mu^{(i)}$ for all $i$. 
This relation factors through the action of $\OLA(\FQT)$ to give a relation $\tri$ on ${\mathfrak L}/O_{\ell'}(\FF_q)$. 

Next, we define a relation $\tri$ on $\IBr(SL_{n})$
as follows. By Corollary~\ref{CLabel},  $\IBr(SL_n)=\{L(\sfr)_j\mid \sfr\in {\mathfrak L}/O_{\ell'}(\FF_q),1\leq j\leq \kappa_{\ell'}({\mathfrak s})\kappa_{\ell}({\mathfrak s})\}$. 
For each $\sfr\in {\mathfrak L}/O_{\ell'}(\FF_q)$ fix an arbitrary {\it linear} order $\tri$ on the set 
$$\{L(\sfr)_j\mid 1\leq j\leq \kappa_{\ell'}({\mathfrak s})\kappa_{\ell}({\mathfrak s})\}.$$ 
If $\sfr,\tfr\in {\mathfrak L}/O_{\ell'}(\FF_q)$, $\sfr\neq\tfr$, set $L(\sfr)_i \tri L(\tfr)_j$ if and only if  $\sfr \tri \tfr$.

\begin{Lemma}\label{order}
$\tri$ is a partial order on  ${\mathfrak L}/O_{\ell'}(\FF_q)$ and on
$\IBr(SL_{n})$.
\end{Lemma}

\begin{proof}
First, consider $\tri$ on ${\mathfrak L}/O_{\ell'}(\FF_q)$. Clearly, it is 
reflexive and transitive. Assume in the above notation that $\sfr \tri \tfr$ and 
$\tfr \tri \sfr$. Then $a = b$, and, after a suitable relabeling,  
$[\si_{i}\nu] = [\tau_{i}]$ and $\la^{(i)} \tri \mu^{(i)}$ for 
all $i$ and for some $\nu \in \OLA(\FQT)$. Also, there is  $\pi \in \SSS_{a}$ such that
$\mu^{(\pi(i))} \tri \la^{(i)}$ for all $i$. Thus 
$$\la^{(1)} + \ldots + \la^{(a)} \tri \mu^{(1)} + \ldots + \mu^{(a)} = 
  \mu^{(\pi(1))} + \ldots + \mu^{(\pi(a))} \tri \la^{(1)} + \ldots + \la^{(a)}.$$
By Lemma \ref{obv}(ii), $\la^{(i)} = \mu^{(i)}$ for 
all $i$, i.e. $\sfr$ and $\tfr$ belong to the same $\OLA(\FQT)$-orbit. Thus 
$\tri$ is anti-symmetric.
Now, $\tri$ on $\IBr(SL_{n})$  is clearly 
reflexive and transitive. If $L(\sfr)_i \tri L(\tfr)_j$ and 
$L(\tfr)_j \tri L(\sfr)_i$, we may assume $\sfr\neq\tfr$. Then $\sfr \tri \tfr$ and $\tfr \tri \sfr$, and so  $\sfr=\tfr$, a contradiction.
\end{proof}

\begin{Proposition}\label{cf-sl}
Let $G = GL_{n}$, $S = SL_{n}$. For any $V \in \IBr(G)$ there is $\VC \in \Irr(G)$ such that:
\begin{enumerate}
\item[{\rm (i)}] $V$ is a composition factor of $\VC (\mod \ell)$;
\item[{\rm (ii)}] $\kappa^{G}_{S}(V) = \kappa^{G}_{S}(\VC)$;
\item[{\rm (iii)}] If $W$ is an irreducible constituent of $V\dar_{S}$ then there is 
a unique irreducible consitutuent $\WN$ of $(\VC)\dar_{S}$ such that $W$ is a multiplicity  one composition factor of 
$\WN (\mod \ell)$ and all other composition factors $U$ of
$\WN (\mod \ell)$ satisfy $U \tri W$.      
\end{enumerate}
\end{Proposition}
\begin{proof}
Write $V = L(\sfr)$, 
 $\kappa_{\ell}(\sfr) = \ell^{c}$ for $c\in\ZZ_{\geq 0}$, and 
$\kappa_{\ell'}(\sfr) = |J|$, where 
$J := \{ \nu \in \OLA(\FQT) \mid \nu \cdot \sfr = \sfr\} = \langle \rho \rangle$.
As $\ell^{c}$ divides 
each $(\la^{(i)})'$, we may choose $\mu^{(i)}$ with
$(\la^{(i)})' = \ell^{c}(\mu^{(i)})'$. Next, consider a  
$\rho$-orbit 
\begin{equation}\label{frho}
  ([\si_{j_{1}}],\la^{(j_{1})}) \mapsto ([\si_{j_{2}}],\la^{(j_{2})}) \mapsto \ldots 
  \mapsto  ([\si_{j_{b}}],\la^{(j_{b})}) \mapsto ([\si_{j_{1}}],\la^{(j_{1})}).
\end{equation}
By Lemma \ref{parts}, the $\si_{j_{i}}$ all have the same degree, say $d$. Also, because of 
(\ref{frho}) we may assume that $\si_{j_{i}} = \si_{j_{1}}\rho^{i-1}$ for $1 \leq i \leq b$. 
Let $u: = u(c,d)$ be as in Lemma~\ref{k2c}. Set
$  \tau_{j_{m}} = \si_{j_{m}}u(c,d)$ for $1\leq m\leq b$. 
Do the same for all $\rho$-orbits on the
set of components of $\sfr$, and set $\VC := \SC(\tfr)$, where
$$\tfr := [([\tau_{1}],\mu^{(1)}),([\tau_{2}],\mu^{(2)}), \ldots ,([\tau_{a}],\mu^{(a)})].$$ 
As $(\tau_{i})_{\ell'}=\si_{i}$ and 
$\deg(\tau_{i}) = \ell^{c}\deg(\si_{i})$ for all $i$, (i) holds by Theorem~\ref{cf-main}.

Next we show that $\kappa_{\ell'}(\VC) = \kappa_{\ell'}(\sfr)$. If $I:=\{ \nu \in \OLA(\FQT) \mid \nu \cdot \tfr = \tfr\}$ then $\kappa_{\ell'}(\VC)=|I|$. If   
$\nu \in I$  sends 
$([\tau_{i}],\mu^{(i)})$ to $([\tau_{j}],\mu^{(j)})$ then $\mu^{(i)} = \mu^{(j)}$ and 
$[\nu\si_{i}] = [\si_{j}]$ by Lemma \ref{parts}(ii). Hence 
$\la^{(i)} = \la^{(j)}$ 
 $\nu$ sends $([\si_{i}],\la^{(i)})$ to $([\si_{j}],\la^{(j)})$, i.e. 
$\nu \in J$. Thus $I\subset J$. To get the opposite inclusion we show that $\rho\in I$. This is obvious if 
$\ell {\not{|}} (q-1)$ as then $\tfr = \sfr$. 
So  assume $\ell|(q-1)$ and consider a $\rho$-orbit as in (\ref{frho}). Then $\la^{(j_{1})} = \ldots = \la^{(j_{b})}$, and so
$\mu^{(j_{1})} = \ldots = \mu^{(j_{b})}$. Next, it is clear that 
$\rho$ sends $\tau_{j_{i}}$ to $\tau_{j_{i+1}}$ for $1 \leq i \leq b-1$. It remains to show
$[\si_{j_{1}}\rho^{b}u(c,d)] = [\si_{j_{1}}u(c,d)]$. Denote 
$\si := \si_{j_{1}}$, $\beta := \rho^{b}$, $u := u(c,d)$, $k:=|u|$. Then 
$[\si \beta] = [\si]$ by (\ref{frho}), whence $\si\beta = \si^{q^{e}}$ and 
so $\beta = \si^{q^{e}-1}$ for some $e\in \ZZ_{\geq 0}$. Therefore
$$\beta^{k} = \beta^{1 + q^{e} + \ldots + q^{(k-1)e}} = 
  \si^{q^{ke}-1} = (\si u)^{q^{ke}-1}.$$
The last equality holds because $\ell|(q-1)$ and so $k$ divides $q^{ke}-1$ by Lemma \ref{LNT}.
Thus $[\si u] = [(\si u)^{q^{ke}}] = [\si u\beta^{k}]$, i.e. $\beta^{k}$ preserves 
$[\si u]$. But $\beta$ is an $\ell'$-element and $k$ is an $\ell$-power, hence $\beta$ also 
preserves $[\si u]$, as required.
 
Now we show that $\kappa_{\ell}(\VC) = \kappa_{\ell}(\sfr) $. Recall that
$[(\tau_{i})_{\ell'}] = [\si_{i}] \neq [\si_{j}] = [(\tau_{j})_{\ell'}]$ whenever $i \neq j$. 
So by Lemma~\ref{parts}(i),  
$\gamma \in \OLB(\FQT)$ preserves 
$\tfr$ precisely when $\gamma$ preserves each $[\tau_{i}]$, which by Lemma \ref{k2c} happens exactly when $|\gamma|$ divides 
$\ell^{c}$, whence the claim. In particular, we have established (ii). 

By \cite{DJ2}, in the Grothendieck group of $\FF G$-modules, we have 
$$\overline{\SC(\tau_{i},\mu^{(i)})} = L(\tau_{i},\mu^{(i)}) 
  + \sum_{\nu^{(i)} \tri \mu^{(i)},~\nu^{(i)} \neq \mu^{(i)}}
    x^{i}_{\nu^{(i)}}L(\tau_{i},\nu^{(i)})$$
for $x^{i}_{\nu^{(i)}} \in \ZZ$. As 
$[(\tau_{i})_{\ell'}]=[\si_{i}]$ and
$\deg(\tau_{i}) = \ell^{c}\deg(\si_{i})$, by \cite{DJ1} we have
$$\overline{\SC(\tau_{i},\mu^{(i)})} = L(\si_{i},\la^{(i)}) 
  + \sum_{\nu^{(i)} \tri \mu^{(i)},~\nu^{(i)} \neq \mu^{(i)}}
    x^{i}_{\nu^{(i)}}L(\si_{i},\eta^{(i)}),$$
where $(\eta^{(i)})' = \ell^{c}(\nu^{(i)})'$. As $\nu^{(i)} \tri \mu^{(i)}$, we have 
$(\la^{(i)})' = \ell^{c}(\mu^{(i)})' \tri \ell^{c}(\nu^{(i)})' = (\eta^{(i)})'$, 
i.e. 
$\eta^{(i)} \tri \la^{(i)}$. As everything is true for every $1\leq i\leq a$, we have
\begin{equation}\label{EFinal}
\bar{V}_{\CC} = V + \sum_{\eta^{(1)}, \ldots ,\eta^{(a)}}y_{\eta^{(1)}, \ldots ,\eta^{(a)}}
  L(\si_{1},\eta^{(1)}) \circ \ldots \circ L(\si_{a},\eta^{(a)}),
\end{equation}
for some $y_{\eta^{(1)}, \ldots ,\eta^{(a)}} \in \ZZ$. Moreover, $\eta^{(i)} \tri \la^{(i)}$ for all $i$ and $\eta^{(j)} \neq \la^{(j)}$ for at 
least one $j$.
Furthermore, each $V' := L(\si_{1},\eta^{(1)}) \circ \ldots \circ L(\si_{a},\eta^{(a)}) \in \IBr(G)$.
By definition of 
$\tri$, for any irreducible constituent $W'$ 
of $V'\dar_{S}$ and for any irreducible constituent $W$ of $V\dar_{S}$, we have $W' \tri W$ and 
$W' \not\simeq W$. 

Also, by (ii), we may 
write $(\VC)\dar_{S} = \sum^{\kappa}_{i=1}W^{i}_{\CC}$ and 
$V\dar_{S} = \sum^{\kappa}_{i=1}W^{i}$, where $W^{i}_{\CC} \in \Irr(S)$, 
$W^{i} \in \IBr(S)$, and $\kappa = \kappa^{G}_{S}(V)$. We may assume that 
$W^{1}$ is a composition factor of $W^{1}_{\CC} (\mod \ell)$. Since 
$W^{1}_{\CC}$ and $W^{1}$ have the same inertia group in $G$, we may relabel the 
$W^{i}_{\CC}$ so that $W^{i}_{\CC} = \ ^{g_{i}}W^{1}_{\CC}$ 
whenever $W^{i} = \ ^{g_{i}}W^{1}$. Hence $W^{i}$ is a composition factor of 
$W^{i}_{\CC} (\mod \ell)$ for all $i$. 
Now (iii) follows from (\ref{EFinal}) and the following remarks.
\end{proof}

\medskip
\noindent
{\bf Proof of Theorem \ref{main3}.} Let $W$ be an irreducible $\FF SL_{n}$-module. 
Then $W$ is a constituent of $V\dar_{SL_{n}}$ for some 
$V \in \IBr(GL_{n})$. Now we apply Proposition \ref{cf-sl} to $W$ and 
define $\Theta(W) = \WN$. 
\hfill $\Box$

\smallskip
Taking intersection with any $\ell$-block $B$ of $SL_{n}$, we also get a block version
of Theorem \ref{main3}: {\it The decomposition submatrix of $B$, with respect to 
$\Theta(B \cap \Irr(SL_{n})$ and the order $\tri$ on $B \cap \IBr(SL_{n})$, is 
lower unitriangular.} In particular, 
$\Theta(B \cap \Irr(SL_{n}))$ yields an ordinary basic set for $B$.

\vspace{2mm}

\small
\ifx\undefined\bysame
\newcommand{\bysame}{\leavevmode\hbox to3em{\hrulefill}\,}
\fi

\end{document}